\documentstyle[12pt]{article}

\newfont{\goth}{eufm6 at 12pt}

\def\vespace{\rule{0pt}{20pt}}

\def\R{{\rm I\hspace{-2.5pt} R}}
\def\C{{\mbox{\rm\bf C\hspace{-7.3pt}}{^{_{\bf\mid}}}\hspace{4.5pt}}}
\def\D{{\mbox{\rm\bf D\hspace{-7.3pt}}{^{_{\bf\mid}}}\hspace{4.5pt}}}
\def\car{\par\hfill$\rule{3mm}{3mm}$\par}

\def\mapright#1{\smash{\mathop{\makebox[12mm]{\rightarrowfill}}\limits^{\scriptstyle{#1}}}}
\def\mapleft#1{\smash{\mathop{\makebox[12mm]{\leftarrowfill}}\limits^{\scriptstyle{#1}}}}
\def\mapdown#1{\Bigg\downarrow \rlap{$\vcenter{\hbox{$\scriptstyle{#1}$}}$}}
\def\mapup#1{\Bigg\uparrow \rlap{$\vcenter{\hbox{$\scriptstyle{#1}$}}$}}
\def\mapvegal#1{\left|\right|  
{\vbox{\vsize=6mm}}\rlap{$\vcenter{\hbox{$\scriptstyle{#1}$}}$}}

\def\div{\mbox{div}}
\def\rdc{{\displaystyle{{\widehat{\cal{R}}}_2}}}
\def\rtc{{\displaystyle{{\widehat{\cal{R}}}_3}}}
\def\fuc{{\displaystyle{{\hat{f}}_1}}}
\def\fdc{{\displaystyle{{\hat{f}}_2}}}
\def\ftc{{\displaystyle{{\hat{f}}_3}}}
\def\ftci#1{{\displaystyle{{\hat{f}}{}^{-1}_#1}}}

\def\Jus{{\displaystyle\underline{J_1}}}
\def\Jds{{\displaystyle\underline{J_2}}}
\def\jtm#1{{\displaystyle\underline{{J_{#1}}(T{\cal M})}}}

\def\rsu{{\displaystyle\not{\!\!{\cal{R}}_{{}_1}}}}

\def\rsd{{\displaystyle\not{\!\!{\cal{R}}_{{}_2}}}}
\def\msd{{\displaystyle\not{\!\!{M}}_{{}_2}}}

\def\gtm#1{{\displaystyle{\Lambda^{#1}}\,{T^\ast}\!\!{\cal M}}}
\def\gtms#1{{\displaystyle\underline{{\Lambda^{#1}}\,{T^\ast}\!\!{\cal M}}}}

\def\otith{{\displaystyle{\otimes_{{}_{\!{\theta}}}}}}

\def\stm#1{{\displaystyle{S_{#1}}{T^\ast}\!\!{\cal M}}}
\def\stms#1{\underline{{S_{#1}}{T^\ast}\!\!{\cal M}}}
\def\tsm{{\displaystyle{{T^\ast}\!\!{\cal M}}}}
\def\tsms{{\displaystyle\underline{{T^\ast}\!\!{\cal M}}}}

\def\twinu{{\displaystyle{\bf D}_1}}
\def\twind{{\displaystyle{\bf D}_2}}
\def\twinup{{\displaystyle{\bf D'}_1}}
\def\twindp{{\displaystyle{\bf D'}_2}}
\def\twinupc{{\displaystyle{\bf \widehat{D}'}_1}}
\def\twindpc{{\displaystyle{\bf \widehat{D}'}_2}}
\def\twinuc{{\displaystyle{\bf \widehat{D}}_1}}
\def\twindc{{\displaystyle{\bf \widehat{D}}_2}}

\def\wedgebar{{\displaystyle\wedge^{\hskip -2.6mm -}}}
\def\twinsu{{\displaystyle\bf\not{\!\!D}'_1}}
\def\itwinsu{{\displaystyle\bf\not{\!d}'_1}}
\def\ditwinsu{{\displaystyle{{}^\ast}\!\!\!\bf\not{\!d}'_1}}
\def\twinsd{{\displaystyle\bf\not{\!\!D}'_2}}
\def\mc#1{{\displaystyle\widehat{M}_{#1}}}
\def\de{{\displaystyle{\bf d}}}

\def\cc#1{{\displaystyle{\widehat{\chi}^{(#1)}}}}
\def\ccc#1{{\displaystyle{\widehat{\chi}_{#1}}}}

\def\tc#1{{\displaystyle{\widehat{\tau}^{(#1)}}}}
\def\tpc#1{{\displaystyle{\tau}^{(#1)}}}
\def\tcc#1{{\displaystyle{\widehat{\tau}_{#1}}}}

\def\pve{{\displaystyle\widehat{\cal{A}}}}
\def\ptge{{\displaystyle\widehat{\cal{B}}}}

\def\grobral{{\displaystyle\,[{{\hskip -1.7mm}\rule[-1.04mm]{1mm}{4.2mm}}\,}}
\def\grobrar{{\displaystyle]{{\hskip -.6mm}\rule[-1.04mm]{1mm}{4.2mm}}\,}}

\def\dgou#1{{\displaystyle\widehat{\mbox{\goth D}}_{#1}}}

\def\as#1#2{{\displaystyle\bf\not{\!\!A}{}^{#1}_{#2}}}
\def\cs#1#2{{\displaystyle\bf\not{{\hskip -1pt}{\cal C}}{}^{#1}_{#2}}}
\def\dcs#1#2{{\displaystyle{{}^\ast}\!\!\!\bf\not{{\hskip -1pt}{\cal  
C}}{}^{#1}_{#2}}}

\begin{document}

\parindent=6mm
\title{\bf From a Relativistic Phenomenology of Anyons to a Model
of Unification of Forces via the Spencer Theory of
Lie Structures}
\author{Jacques L.~Rubin\thanks{Institut du Non-Lin\'eaire de Nice,
UMR 129 C.N.R.S.-Universit\'e de Nice-Sophia-Antipolis,
1361 route des Lucioles, 06560 Valbonne, France.
E-mail: rubin@inln.cnrs.fr. Fax: 00/33/04-92-96-73-33}}
\date{\today }
\maketitle

\begin{abstract} Starting from a relativistic phenomenology of
anyons in crystals, we discuss the concept of relativistic
interaction and the need to unify electromagnetism and gravitation
within the Spencer cohomology of Lie equations. Then, from the
sophisticated non-linear Spencer complex  of the Poincar\'e and
conformal Lie pseudogroups, we build up a non-linear relative
complex assigned to a gauge sequence for
electromagnetic and gravitational potentials and fields. Then,
using a conformally equivariant Lagrangian density, we deduce,
first, the two first steps of its corresponding Janet complex and
second, the dual relative linear complex. We conclude
by giving suggestions for higher unification with the  weak and
strong interactions and interpretations of the Lagrangian density
as a thermodynamical function and quantum wave-function.
\end{abstract}

\vfill\eject
\tableofcontents
\mbox{ }
{\bf keywords:} anyons, superconductivity, unification models,
electromagnetism, gravitation, gauge theory, non-linear Spencer
cohomology,  conformally flat Riemannian manifolds, Lie groupo\"{\i}ds,
G-structures, Lie pseudogroups, deformations, differential
complexes. Lie equations.\\
{\bf Short running title}: From Anyons to Unification of Forces via the
Spencer Theory
\vfill\eject
\addcontentsline{toc}{section}{1. Introduction}
\section*{\centering 1. Introduction}

\hskip 6mm In this paper we propose a model as well as suggestions
for a  unification of physical interactions. This is a model of
electromagnetic and  gravitational interactions well-founded on a
phenomenological relativistic   model of anyons in high-$T_c$
superconductors (Rubin 1994). This unification has its  roots,
first in the Spencer cohomology of the conformal Lie  pseudogroup
which has abundantly been studied by J. Gasqui \& H. Goldschmidt
(Gasqui {\it et al.} 1984),  and second, in the non-linear cohomology of Lie
equations studied by   H. Goldschmidt \& D. Spencer (1976a, 1976b, 1978a,  
1978b, 1981, see references therein). Meanwhile
we only partially refer to
some of its aspects to work out a relative non-linear cohomology
explicitly associated with a model of unification.   Such an approach
was originally proposed by J.-F. Pommaret (1988, 1989, 1994) however, it  
appeared to us to be
incomplete, indeed erroneous.\par  Thus, the purpose of the present
paper is to discuss the Pommaret model and to suggest new
developments based on the same assumption. Like Pommaret
(1989), we think that the geometrical approach of the
Maxwell theory has to be modified to be incorporated in a larger
theory which explicitly includes gravitation through different
equations describing the variations of potentials of gravitation.
This result - or proposition - has not been obtained by Pommaret
(Pommaret 1994, see conclusions page 456) who could not find any
alternative descriptions and justifications neither for the Einstein
theory, nor for the equations of fields of gravitation within the
frame of the Spencer cohomology.\par  We conclude this paper a) by
succinctly proposing  a possible reinterpretation of the quantum
wave-function as a classical thermodynamic function within the
frame of the Misra-Prigogine-Courbage (MPC) ergodic theory of
fields (Misra {\it et al.} 1979, Misra 1987), b) by suggesting ideas about a  
unification
including the weak and strong interactions along a basic ``\`a la
Penrose" approach (Penrose {\it et al.} 1986).\par Furthermore, this work is  
the result of informal reflexions about an increasing amount of
contradictions and incoherences mainly concerning the concept of
relativistic interaction, (that we find more and more serious) in
the field of quantum physics as well as in classical
physics.  According to this observation, we first present our
motives and a description of these contradictions in relation to
F. Lur\c cat's (1964), J.-M. L\'evy-Leblond's (1990)
and J.-F. Pommaret's arguments (1989) in order to justify a
development via the Spencer cohomology of Lie pseudogroups (Goldschmidt {\it  
et al.} 1976a, 1976b, 1978a, 1978b, 1981).

\addcontentsline{toc}{section}{2. Goals and problems}
\section*{\centering 2. Goals and problems}

\addcontentsline{toc}{subsection}{2.1. The physics of crystals and a
relativistic phenomenology of anyons}
\subsection*{\centering 2.1. The physics of crystals and a relativistic
phenomenology of anyons}

\hskip 6mm Our initial motivation shall be seen as extremely far
from the problems with unifications. Actually, we were more concerned
in a simple minor model of a relativistic phenomenology of creation
of anyons, accurate for certain crystals (Rubin 1994).  At the
origin of this process of creation, we suggested the
kinetico-magnetoelectrical effect as described by  E. Asher
(1973) and which has its roots in the former Minkowski
works about the relations between tensors of polarization $P$ and
Faraday tensors $F$ in a moving material of optical index
$n\not =1$.These relations are established by turning the following
diagram  into a commutative one:


$$\begin{array}{rcl}
F'&\mapright{\Lambda}&F\cr
\mapdown{\chi'}&&\mapdown{\chi}\vespace\cr
P'&\mapright{\Lambda}&P\vespace
\end{array}$$


where $\Lambda$ is a Lorentz transformation, allowing us to shift
from a frame
$R'$ to a frame $R$, and $\chi'$ and $\chi$ are respectively  the
tensors of susceptibility within those two frames, as well
supposing $P$ (or $P'$)  linearly depending on $F$ (respectively
$F'$). Resulting from this commutativity, the tensor $\chi$
linearly depends on $\chi'$ in general and also on a velocity
4-vector $\tilde u$ associated to $\Lambda$ (i.e. the relative
velocity 4-vector between $R$ and $R'$). In assimilating $R'$ to
the moving crystal frame and $R$ to the laboratory frame, then to
an applied electromagnetic field $F$ fixed in $R$, corresponds
in $R'$ a field of polarization $P$ which varies in relation to
$\tilde u$. This is the so-called kinetico-magnetoelectrical
effect. \par Parallel to this phenomenon,  A. Janner \& E.
Asher  studied the concept of relativistic point symmetry in
polarized crystals (Janner {\it et al.} 1969, 1978). Such a symmetry is  
defined, on the one hand, by a given discret group $G$, sub-group of the
so-called Shubnikov group $O(3)1'$ associated with the crystal, and
on the other hand, as satisfying the following
properties: to make this relativistic symmetry exist, there must
be a $H(P)$ non-trivial group of Lorentz transformations
depending on $P$, in which $G$ is a normal  sub-group, and that
leaves the tensor of polarization $P$ invariant. In other words, if
$N(G)$ is the normalizer of $G$ in the Lorentz group $O(1,3)$, and
$K(P)$ the sub-group of $O(1,3)$ leaving $P$ invariant, then $H(P)$
is the maximal sub-group such that:
	

$$\left\{ \begin{array}{ll}
H(P)\subseteq K(P)\cap N(G)\\
H(P)\cap O(3)1'=G.
\end{array}\right.$$


We can prove that $H(P)$ is about to exist only if a particular
non-vanishing set $V$ of velocity 4-vectors, invariant by action of
$G$, is present and consequently compatible with a
kinetico-magnetoelectrical effect (Asher 1973). Therefore, if there
is an interaction between moving particles in the crystal and the
polarization $P$, then the trajectories and $P$ are obviously
modified, and so is $H(P)$. In this process, only the group $N(G)$
is conserved so that the polarization and the trajectories are
deducible  during the time by the action of $N(G)$.\par  As we
shall stipulate later on, the existence of an interaction will
emerge due to a correlation between the  position 3-vectors $\vec
r$ of the charge carriers and a particular 3~-~vector $\vec w$
($\notin V$ in general) associated with $P$; $\vec w$ becoming then a
function of $\vec r$. In order to allow a cyclotron-type motion
which is implicit within the theory of anyons, the group $N(G)$
must contain the group $SO(2)$ and the latter must non-trivially
act on all the groups $H(P)$ associated to $G$. Then, only 12
groups $G$ are compatible with such a description (Rubin 1993). In
fact, throughout this development, we implicitly use a principle of
equivalence similar to the one formulated in general relativity:
one cannot distinguish a cyclotron~-~type motion in a constant
polarization field from a uniform rectilinear motion in a field of
polarization varying in time by action of the normalizer $N(G)$.
The time evolution of $\vec w$ requiring the explicit knowledge of
the gradient ${\vec \nabla}({\vec w})$ of $\vec w$, it follows that
$\vec w$ and ${\vec\nabla}({\vec w)}$ are respectively the analogues
(not the equivalents) of the tetrads and of the Killing vector
fields. From an other point of view, the interaction is
considered to allow the extension of an invariance with respect
to $H(P)$ to an invariance with respect to $N(G)$. The lack of
interaction is then what breaks down the symmetry!\par This type
of reasoning concerns in fact a large amount of physical
phenomena such as the spin-orbit interaction for instance. In
this context, the cyclotron-type motion of electrons in anyonic
states would be similar to the Thomas or Larmor precessions (see
also the Coriolis or Einstein-Bass effects). More precisely,
taking up again a computation, analogous to the Thomas
precession one (Bacry 1967) (i.e.  considering as a constant the
scalar product of two tangent vectors being two parallel
transports along the trajectory (Dieudonn\'e 1971)), concerning a
charge carrier with the velocity  4-vector $\tilde u$  in $R$,
``polarized'' by
${\vec w}({\vec r})$ such as for example  (${\tilde v}=(0,{\vec
v})_R$ constant and ${\vec v}\in V$):


$$
{\tilde w}=(0,{\vec w})_R\equiv -P.{\tilde v}\ or\ {}^{\ast}P.{\tilde v},
$$


where $P$ depends on $\vec r$, one can prove from ${\tilde
w}.{\tilde u}=cst.$ that ($t$ being the laboratory frame time
and ($\tilde{r}=(t,\vec{r})_R$):


\begin{equation}
{{d\tilde{u}}\over{dt}}= (-e/m)F_{eff.}(\tilde{r}).\tilde{u},
\label{i}
\end{equation}


where $m$ and $e$ are respectively the mass and the electric charge
of the carrier and $F_{eff.}(\tilde{r})\equiv ({\vec E}_{eff.}
(\tilde{r}),{\vec B}_{eff.}(\tilde{r}))$ is an effective
Faraday tensor such that ($\gamma = (1-{\vec w}^2)^{-{1\over 2}}$ and
${\vec j} =e{\vec u}$):

$$
\left\{
\begin{array}{ll}
{\vec B}_{eff.}(\tilde{r})=(m/e^2)\left({\gamma\over 1+\gamma}\right){\vec w}
\wedge\left[{\vec j}.{\vec \nabla}\right]{\vec w}\\
{\vec E}_{eff.}(\tilde{r})={\vec 0}.
\end{array}
\right.
$$


Clearly, $F_{eff.}$  is an element of the Lie algebra of the group $SO(2)$
included in $N(G)$ and with ${\vec B}_{eff.}\in V$. Therefore this ${\vec
B}_{eff.}$ magnetic field or $\vec v$ (up to a constant) might be considered
as the effective magnetic field of the flux-tube $V$ generating
the so-called Aharanov-B\"ohm effect at the origin of the
statistical parameter in the anyons theory (Wilczek 1990). In this
precession computation, from a more mathematical point of view,
taking up the Lie groupo\"{\i}ds theory, we think that perhaps
we shift from a source ``description'' (at $t=0$) to a target one
(at any $t\not =0$). It is definitely an equivalence principle
analogous to the one occuring in general relativity, as shown by
J.-F. Pommaret (1989). Let us add that in general
$\div({\vec B}_{eff.})\not =0$ so that one gets a non-vanishing
density of effective magnetic monopoles generated by the local
variations (due to the interaction) of the polarization vector
field ${\vec w}$ in the crystal. Thus an anyon would be an
effective magnetic monopole associated with a charge  carrier,
namely a dyon. Moreover, because this effective Faraday tensor
is no more a closed two-from, a non-vanishing Chern-Simon has to
be taken into account in a Lagrangian description of anyons,
from which non-vanishing spontaneous constant currents can
occur.

\addcontentsline{toc}{subsection}{2.2. Relativistic interaction
in quantum mechanics}
\subsection*{\centering 2.2 Relativistic interaction in quantum mechanics}

\addcontentsline{toc}{subsubsection}{2.2.a. Polarization in
quantum mechanics}
\subsubsection*{\centering 2.2.a. Polarization in quantum mechanics}

\hskip 6mm The previous vector  ${\vec B}_{eff.}$ (or ${\vec v}\in
V$) has the same status as the spin.  Like the latter, it is
defined by a torsor of  order 2. From that time on, the transition
from a classical description to  a quantum one means that one must
give an account of the interaction  between a free particle
(constituting a first sub-system) and a torsor field of order 2
(constituting a second sub-system). The problem seems to be solved
and in particular the spin appears  to be a ``minor" complication
of the wave-function defined on the Minkowski  space-time, i.e. on
a ``non-polarized" space. This has to be taken as a  postulate, an
erroneous one according to part of F. Lur\c{c}at's (1964) and
J.-M.  L\'{e}vy-Leblond's
(1990) arguments.\par  In fact, it is not even the
case according to classical Galilean mechanics. Considering a
body at a given time, one needs to know 6 parameters to describe
it: 3 for the position and 3 others for the orientation. The
latter are forgotten during the transition from the ``extended
body to the punctual  particle" according to the quantum
description. This can only be justified providing that the energy
of rotation is negligible compared to the energy of translation.
This is the indicator of an inadequacy of the principle of
correspondence from classical mechanics to quantum mechanics. But
the transition from quantum to classical mechanics is
equally  problematic: reaching the limit $\hbar\rightarrow 0$, the
spin vanishes.  Finally, we must say that this correspondance does
not exist any more in chromodynamics.\par Therefore, the fact that
the wave-function only  depends on the position should rather be
considered as a postulate (moreover,  the fact is inexistant in
classical mechanics). For instance, if we measure  the electric and
magnetic fields at a point in space-time with a system of
coordinates defined by a given Lorentz transformation, we can
deduce that  the Faraday tensor is a function of the 4 parameters
of the position and  6 others defining $\Lambda$. It explicitly
occurs within the  kinetico-magnetoelectric effect.\par  Refering now
to the wave-equations and to the methods usually accepted to
determine them, we then use a principle of invariance. First, we
identify  the appropriate space-time symmetries, that is the
group of relativity of the theory (for example the Poincar\'{e}
group). Then, according to Wigner's theory, we build the
irreducible unitary representations  (eventually projective) of the
group to which correspond the elementary ``kinetic" objects of the
theory (i.e. the vector-valued wave-functions). Then again, we
derive the infinitesimal generators of the group that we identify with
physical and geometrical observables like the energy, the kinetic
moment, etc...; the Lie algebra of the group and its  ``quantum"
extensions defining the commutators.  Finally, determining the
invariants through the action of the group we obtain the other
physical observables and their commutators.\par  Unfortunately, some
ambiguities appear. For example, the generalization of the Dirac
theory for spins higher than 1/2 gives different non-equivalent
possible wave-equations. Lastly, if we consider a theory of
interacting fields, the particles associated with these fields can
get out  of their mass shell but not out of their spin shell!  We then
forget the spin again, which is impossible within a ``(m,s)"
theory. Therefore some of the aspects of the interaction and of
the wave-equations - to be brief- can fully account for
neither the free particles, nor the interacting particles!\par  In
order to escape from these contradictions and to best describe the
free particles, F.\,~Lur\c{c}at proposed an approach ``\`{a} la
Wigner". In this way, he postulated that first the scalar complex
wave-function  of a free particle was defined on the ``Poincar\'{e}
space" of the Poincar\'{e}  Lie group. Second, this wave-function
is an eigenfunction of the two Casimirs of  the Poincar\'{e} group.
Thus, the wave-function $\phi$ is a function first, of a 4-vector
position $\tilde{x}$, element of the dual Lie algebra of the  group
of translations and second, of a second order torsor $F$ (like the
Faraday torsor which is an element of the dual Lie algebra of the
Lorentz group).  Also in fact, the Poincar\'{e} Lie group on the
``Poincar\'{e} space" is the action  of a so-called Lie
groupo\"{\i}d on its associated Lie algebro\"{\i}d, defined on
the Minkowski space-time. Finally, if we want to describe the
interaction,  making these modifications and requiring a gauge
invariance, the  involution of the infinitesimal generators
``deformed" by the gauge fields  is no longer satisfied, nor are
the relativistic invariance and the  correspondence principle.\par
This breaking of the invariance can be seen with the Dirac
equation. More exactly, the Dirac equation is not {\it
equivariant},  but only covariant, meaning that this equation is not
invariant  under any conformal changes of coordinates, but only
under a change of frame, i.e. a change of basis of the tangent
Minkowski space. In that case, for instance, one of the spectacular
manifestations  is Klein's paradox of non-conservation of the
current of probability during  the scattering in a square potential
(Itzykson {\it et al.} 1980). Still, the idea stays that $\phi$ is  defined on a
``larger" space-time than the Minkowski one, but for being  a
Kaluza-Klein type theory for example.

\addcontentsline{toc}{subsubsection}{2.2.b. Concerning the contradictions of an 
approach ``\`{a} la Wigner" and the Einstein-Cartan unification}
\subsubsection*{\centering 2.2.b. Concerning the contradictions of an
approach ``\`{a} la Wigner" and the Einstein-Cartan unification}

\hskip 6mm Let us consider a complex wave-function $\phi$ depending
on $\tilde{x}$ and $F$ on which the Poincar\'{e} group projectively
acts. Because, then, of the phase arbitrariness, the ten
infinitesimal generators of the group are defined with ten
arbitrary gauge potentials, so-called  Poincar\'{e} potentials.  In
order to keep the Lie algebra structure, i.e. the involution of the
deformed  generators, the fields associated with these potentials
must satisfy some constraint equations. The electromagnetic field
associated with the translations of the group has especially to be
vanishing, which is first absolutly contradictory with $\phi$ as a
function of $F$. Second, that would mean that one cannot describe
any interaction with a field $F$ without a Poincar\'{e} symmetry
breaking. Nevertheless this symmetry is necessary to keep the
relativistic equivariance.\par Moreover, $\phi$ is independent of
the gauge potentials although $\phi$ should be parametrized by
these functions refering to the formal theory of system of partial
differential equations (PDE). In fact, according to the
K\"{a}hler-Cartan theorem, the analytical solutions of an
involutive system of PDE (so formally integrable) depend on certain
constants of integration (Dieudonn\'e 1971), but, contrary to the ODE
(ordinary differential equations),  they depend on arbitrary
$C^\infty$ functions only constrained to verify  Cauchy initial data
(see also Shih  1986, 1987, 1991). These functions should be identified with  
the gauge potentials and with their fields. Still, we can say that the
incoherences are not over!\par Let us assume the presence of a
relativistic interaction, then there consequently exists, if we refer
to the equivalence principle, a proper frame in which $\phi$ is
stationnary. If $\tau$ is the ``proper" time then
${\partial_\tau}\phi=0$.  Within the ``laboratory" frame, we will
get a different equation, moreover, there won't be any relativistic
equivariance; a problem that is similarly encountered with either the
Newton-Wigner position operator or the Dirac position operator. Indeed
we will still be unable to write down this equation (!) for the
reason that to shift to the Laboratory frame, one needs to know the
classical motion of the particle, i.e. to know the equation that
determines the evolution of the 4-vector speed $\tilde u$ identified
with the basic time type 4-vector ${\tilde{e}}_0$ of the proper frame.\par   
Also at this point, the
equation should be established in a certain system of coordinates
and from ${\tilde{e}}_0$ only, one should obtain the other basic
vectors of space type ${\tilde{e}}_i$ (i=1,2,3) as well as their
time evolutions. This equation should then be  integrable in the
Fr\"{o}benius sense. We think it is perhaps a matter of a
generalisation of the Fr\^{e}net moving frame method, as formulated
by E. Cartan.  Whereas we observe that by determining this moving
frame from
${\tilde{e}}_0$  and ${\dot{\tilde{e}}}_0$, that the
${\dot{\tilde{e}}}_i$'s are defined from the ${\tilde{e}}_\mu$'s
($\mu$=0,1,2,3) and from the third order time derivatives  of
${\tilde{e}}_0$. Then, if ${\dot{\tilde{e}}}_0\equiv
F.{\tilde{e}}_0$, we need to know the time
derivatives of $F$ up to the second order. Therefore, in the
laboratory frame, $\phi$
should depend on the derivatives of $F$, contrary  to the
initial assumption, unless these derivatives are themselves
functions of $F$  and $\tilde{x}$. Condition which, first, is not
the case in the Maxwell theory, second would set down some
constraints on the moving frame and consequently a ``partial"
equivariance; unless one completes it with other fields - like those
associated for example with the space-time curvature  -  which would
be related to the  derivatives of $F$.
Let us remark also that this kind of discussion seems to be very
similar to the one encountered in the demonstration of the
Cauchy-Kowalewska theorem when starting from a pfaffian system to a
``normal system" of PDE's and considering the so-called ``regular"
tangent spaces and integrable manifolds (Dieudonn\'e 1971).
\par On that subject, one can notice that the equation (\ref{i}) can
be rewritten  in an orthonormal system of local coordinates
($\alpha, \beta, \gamma= 1, ..., n$; the $\Gamma$'s being the
Christoffel symbols):


$${{\dot u}^\alpha}+
{\Gamma{}_{\beta,\gamma}^{\alpha}}{u^\beta}{u^\gamma}=0.$$


We recognize the equation of the geodesics associated to a
Riemannian connexion with torsion. This would suggest a unification
in reference to the Einstein-Cartan theory. If we then keep on with
the assumption that one has to add to the electromagnetic field a
gravitational field and that the derivatives of the fields are
functions of the fields themselves (as with the Bianchi identities
according to the non-abelian theory for example), that means we
make the assumption of the existence of a differential sequence. In
electromagnetism, it is a matter of the de Rham sequence but
gravitation does not interfere. The sequence integrating the latter
- and being the purpose of this paper -  might be a certain
generalizing complex like the Spencer one,  following then a method
proposed by J.-F. Pommaret (1994 )but largely modified.

\addcontentsline{toc}{section}{3. The Lie conformal pseudogroup associated to  
the unification model}
\section*{\centering 3. The Lie conformal pseudogroup associated to the  
unification model}

\hskip 6mm First of all, let us assume that the group of relativity is not
the  Poincar\'{e} group anymore but the conformal Lie group (we
know from Bateman and Cunningham studies (1910) that  it is
the group of invariance of the Maxwell equations). In particular,
this involves  that no changes occur shifting from a given frame to
a uniformly accelerated relative one. From a historical point of
view, that happened to be the starting point of the Weyl theory
which was finally in contradiction with experimental data and for
various other reasons presented, for
instance by J.-F. Pommaret  in the framework of the Janet and Spencer
complexes (Pommaret 1989). Then, starting from this mathematical
framework, J.-F. Pommaret considered in trying this unification, the
linear Spencer complex defined from the system of finite Lie
equations associated to the Lie pseudogroup of conformal
isometries. Unfortunatly, the system of PDE proposed by
J.-F.~Pommaret is incomplete  and its conditions of use are not
really given. On the other hand, he claimed  the Spencer complex of
the conformal Lie pseudogroup would be the ``unification complex"
(Pommaret 1988),  whereas we merely prove that it would rather be a
relative complex deduced  from an abelian
extension (Godschmidt 1976a, 1976b, 1978a, 1978b, 1981), when reaching the  
conformal Lie  pseudogroup
starting from the Poincar\'{e} one.  Before  tackling these
various complexes, we present and recall a few relations concerning
 the conformal Lie pseudogroup action on some tensors such as the
metric and the  Riemann and Weyl curvatures. Let us first call
$\cal M$, the base space (or space-time),  assumed to be of class
$C^\infty$, of dimension $n\geq 4$, connected, paracompact,  without
boundaries, oriented and endowed with a metric 2-form
$\omega$, symmetric, of class $C^2$ on $\cal M$ and
non-degenerated. We also  assume $\cal M$ to be conformally flat.

\addcontentsline{toc}{subsection}{3.1. The conformal finite Lie equations}
\subsection*{\centering 3.1. The conformal finite Lie equations}

\hskip 6mm These equations are deduced from the conformal action on
the metric. Let us consider ${\hat f}\in{Diff{}_{loc.}^1}(\cal{M})$
and any $\alpha\in{C^0}(\cal{M},{\R})$, then if
${\hat f}\in{\Gamma_{\widehat{G}}}$ ($\Gamma_{\widehat{G}}$ being
the pseudogroup of local  conformal bidifferential maps on
$\cal{M}$), ${\hat f}$ is a solution of the following  system of
PDE (other PDE must be satisfied to completly define
$\Gamma_{\widehat{G}}$):

\begin{equation}
      \left\{
         \begin{array}{ll}
{{\hat f}^\ast}\omega=e^{2\alpha}\omega\\
{\rm and\ with}\ \det(J({\hat f}))\not=0,
         \end{array}
      \right.
\label{1}
\end{equation}

where $J({\hat f})$ is the Jacobian of $\hat f$, and ${\hat
f}^\ast$ is the pull-back of $\hat f$. We denote $\tilde{\omega}$
the metric on $\cal{M}$ such as:

$$
\tilde{\omega}\buildrel{\ def.} \over{\equiv} {e^{2\alpha}}\omega,
$$

and we agree to put a tilde on each tensor or geometrical
``object" relative to or deduced from this metric $\tilde{\omega}$.
Let us notice that the latter depends on a fixed given element
${{\hat f}^\ast}\in{\Gamma_{\widehat{G}}}$. Also, in order to
properly recall the last point, we shall sometimes use an
alternative notation such as:

$$
\tilde{\omega}\buildrel{\ def.}\over{\equiv} {{}^{\hat f}\omega}.
$$

This convention of notation will also be used  on each geomerical
object relative to this metric. Now, doing a first prolongation of
the system (\ref{1}), we deduce other second order PDE connecting
the affine connexion 1-forms of Levi-Civita $\nabla$ and
$\widetilde{\nabla}$ respectively associated to $\omega$ and
$\tilde{\omega}$. To obtain these PDE, we merely start from the
following definition of $\widetilde{\nabla}$: let $X$, $Y$ and $Z$
be any vector fields in
${C^1}(T\cal{M})$, ${\hat f}\in{Diff_{loc.}^2}(\cal{M})$ and
$\alpha~\in~{C^1}~(\cal{M}~,~{\R}~)$, then by definition we have:

\begin{equation}
\begin{array}{rcl}
{\tilde{\omega}}({{\tilde{\nabla}}_X}Y,Z)&=&
{1\over{2}}\{
{\tilde{\omega}}([X,Y],Z)+
{\tilde{\omega}}([Z,X],Y)+
{\tilde{\omega}}([Z,Y],Y)\cr
&&+X.{\tilde{\omega}}(Y,Z)+
Y.{\tilde{\omega}}(X,Z)-
Z.{\tilde{\omega}}(X,Y)\},
\end{array}
\label{2}
\end{equation}

from which we deduce with the relation (\ref{1}) $\forall\, X,Y \in
{C^1}(T\cal{M})$ and
$\forall {\hat f}\in {Diff_{loc.}^2}(\cal{M})$,

\begin{equation}
{{\widetilde{\nabla}}_X}Y={\nabla_X}Y+\de\alpha(X)Y+\de\alpha(Y)X-
\omega(X,Y)\,{{}_\ast}\de\alpha,
\label{3}
\end{equation}

where $\de$ is the exterior differential and ${{}_\ast}\de\alpha$
is the dual vector field of the 1-form $\de\alpha$ with respect to
the metric $\omega$, i.e. such as $\forall\, X\in T\cal{M}$

\begin{equation}
\omega(X,{{}_\ast}\de\alpha)=\de\alpha(X)=<\de\alpha|X>.
\label{4}
\end{equation}

Let us also agree to denote in the sequel $\forall\,
\tilde{x}\in\cal{M}$,
$\forall\,\nu\in{\bigwedge^r}{T^\ast}\cal{M}$ and
$\forall\,{\xi_i}\in{T\cal{M}}$
$(i=1,...,r)$:

\begin{eqnarray*}
{\nu_{\tilde{x}}}({\xi_{1,\tilde{x}}},...,{\xi_{r,\tilde{x}}})=<
\nu(\tilde{x})|{\xi_{1,\tilde{x}}}\otimes...\otimes{\xi_{r,\tilde{x}}}>.
\end{eqnarray*}

Prolonging again and using the definition of the Riemann tensor
$\rho$ associated to $\omega$, i.e. $\forall\,X,Y\in{C^2}(T{\cal
M})$

$$
\rho(X,Y)={\nabla_X}{\nabla_Y}-{\nabla_Y}{\nabla_X}-{\nabla_{[X,Y]}},
$$

one obtains the following relation
$\forall\,X,Y,Z\in{C^2}(T\cal{M})$,
$\forall\,\alpha\in{C^2}(\cal{M},{\R})$ and $\forall\, {\hat
f}\in~{Diff_{loc.}^3}(\cal{M})$,

\begin{equation}
\begin{array}{rcl}
\tilde{\rho}(X,Y).Z&=&\rho(X,Y).Z+
\omega(X,Z){\nabla_Y}({{}_\ast}\de\alpha)-
\omega(Y,Z){\nabla_X}({{}_\ast}\de\alpha)\cr
&&+\left\{
\omega({\nabla_X}({{}_\ast}\de\alpha),Z)+
\omega(X,Z)\de\alpha({{}_\ast}\de\alpha)
\right\}Y\cr
&&-\left\{
\omega({\nabla_Y}({{}_\ast}\de\alpha),Z)+
\omega(Y,Z)\de\alpha({{}_\ast}\de\alpha)
\right\}X\cr
&&+\left\{
\de\alpha(X)\omega(Y,Z)-
\de\alpha(Y)\omega(X,Z)
\right\}{{}_\ast}\de\alpha\cr
&&+\left\{\de\alpha(Y)X-\de\alpha(X)Y\right\}\de\alpha(Z)
\end{array}
\label{5}
\end{equation}

Assuming $\cal{M}$ to be comformally flat, the Weyl tensor $\tau$
associated with
$\omega$ vanishes. Hence, the Riemann tensor $\rho$ can be
rewritten
$\forall\,X,Y,Z,U\in{C^2}(T\cal{M})$ as:

\begin{equation}
\begin{array}{rcl}
\omega(U,\rho(X,Y).Z)&=&{1\over{(n-2)}}
\left\{\omega(X,U)\sigma(Y,Z)-\omega(Y,U)\sigma(X,Z)\right.\cr
&&+\left.\omega(Y,Z)\sigma(X,U)-\omega(X,Z)\sigma(Y,U)\right\},
\end{array}
\label{6}
\end{equation}

where $\sigma$ is the so-called Schouten tensor (Gasqui {/it et al.} 1984)
$\forall\,X,Y\in{C^2}(T\cal{M})$,

\begin{equation}
\sigma(X,Y)={\rho_{ic}}(X,Y)-{{\rho_s}\over{2(n-1)}}\omega(X,Y),
\label{7}
\end{equation}

where $\rho_{ic}$ is the Ricci tensor and ${\rho_s}$ is the Riemann
scalar curvature. Thus, we might consider the relation (\ref{6}) as
the existence of a short exact sequence ``symbolically" written as
``$0\rightarrow\sigma\rightarrow\rho\rightarrow\tau\rightarrow0$"
and perhaps related to a sequence of cohomology spaces of symbols
(Gasqui {\it et al.} 1984, Pommaret 1994).
Consequently, the system of PDE (\ref{5}) can be rewritten as a
first order system of PDE concerning $\sigma$. To do this, we first
define two suitable trace operators,  used in the sequel to obtain
the $\tilde{\rho}_{ic}$ and
$\tilde{\rho}_{s}$ tensors and finally the $\tilde{\sigma}$ tensor.
Let us denote $Tr^1$  the trace operator defined such that for any
vector bundle $E$ over $\cal{M}$ we have:

$${Tr^1}:T{\cal M}\otimes\tsm\otimes{E}\longrightarrow{E},$$

with

$${Tr^1}(X\otimes\alpha\otimes\mu)=\alpha(X)\mu$$

for any $X\in T\cal{M}$, $\alpha\in\tsm$ and $\mu\in E$.
Then, the second trace operator is the natural trace $Tr_\omega$ associated to
$\omega$ and defined by:

$${Tr_\omega}:{\buildrel{2}\over\otimes}\tsm\longrightarrow{\R},$$

such that

$${Tr_\omega}(u\otimes v)=v({{}_\ast}u).$$

Finally, with ${Tr^1}\tilde{\rho}={\tilde{\rho}_{ic}}$ and
${Tr_\omega}{\tilde{\rho}_{ic}}={\tilde{\rho}_s}$, we deduce first from the
relations (\ref{5}) and (\ref{7}) $\forall\,{\hat f}\in {Diff_{loc.}^3}(\cal{M})$, 
$\forall\,X,Y\in {C^2}(T\cal{M})$ and $\forall\,\alpha\in{C^2}(\cal{M},{\R})$,

\begin{equation}
\begin{array}{rcl}
\tilde{\sigma}(X,Y)&=&
\sigma(X,Y)+ (n-2)\left(\de\alpha(X)\de\alpha(Y)
-\omega({\nabla_X}({{}_\ast}\de\alpha),Y)\right.\cr
&&\left.-{1\over{2}}\omega(X,Y)\de\alpha({{}_\ast}
\de\alpha)\right).
\end{array}
\label{8}
\end{equation}

In fact this expression can be symmetrized. Indeed, from the following property
satisfied by $\nabla$:

$$
\omega({\nabla_X}({{}_\ast}\de\alpha),Y)+
\omega({{}_\ast}\de\alpha,{\nabla_X}Y)=
X.\,\omega({{}_\ast}\de\alpha,Y),
$$

and the definition of ${{}_\ast}\de\alpha$, one obtains:

$$
\omega({\nabla_X}({{}_\ast}\de\alpha),Y)=X.\,\de\alpha(Y)-\de\alpha({\nabla_X}Y).
$$

But, from the torsion free property of $\nabla$ and from the relation

$$
\de\alpha([X,Y])=X.\,\de\alpha(Y)-Y.\,\de\alpha(X),
$$

one deduces:

$$
\omega({\nabla_X}({{}_\ast}\de\alpha),Y)=Y.\de\alpha(X)-\de\alpha({\nabla_Y}X).
$$

Now, $\forall\,\alpha\in{C^2}(\cal{M},{\R})$, $\forall\,X,Y\in{C^1}(T\cal{M})$, 
defining $\mu\in{C^0}({S_2}\tsm)$ by:

\begin{equation}
\mu(X,Y)={1\over{2}}\left[X.\,\de\alpha(Y)+Y.\,\de\alpha(X)\right],
\label{9}
\end{equation}

one has the relation:

$$
\omega({\nabla_X}({{}_\ast}\de\alpha),Y)=
\mu(X,Y)-{1\over{2}}\de\alpha({\nabla_X}Y+{\nabla_Y}X).
$$

Then, we can rewrite the first order PDE (\ref{8}),
$\forall\,{\hat f}\in{Diff_{loc.}^3}(\cal{M})$,  
$\forall\,X,Y\in{C^2}(T\cal{M})$ and
$\forall\,\alpha\in{C^2}(\cal{M},{\R})$ as:

\begin{equation}
\begin{array}{rcl}
\tilde{\sigma}(X,Y)&=&\sigma(X,Y)+(n-2)\left(\de\alpha(X)
\de\alpha(Y)\right.\cr
&&-\mu(X,Y)+{1\over{2}}\de\alpha({\nabla_X}Y+{\nabla_Y}X)\cr
&&\left.-{1\over{2}}\omega(X,Y)\de\alpha({{}_\ast}\de\alpha)\right).
\end{array}
\label{10}
\end{equation}

It is worthy of note that if ${\rho_s}={c_0}\ ({c_0}\in {\R})$, then
the Schouten tensor $\sigma$ satisfies the relation:

\begin{equation}
\sigma={c_0}{(n-2)\over{2}}\omega.
\label{11}
\end{equation}

and only in this case, the equation (\ref{10}) must become the equation
\ref{1} ans so it desappears.
Then, considering the system (\ref{1}), the system (\ref{10}) reduces to
a second order system of PDE such as $\forall\,\alpha\in{C^2}(\cal{M},{\R})$ and 
$\forall\,X,Y\in{C^1}(T\cal{M})$ we have:

\begin{equation}
\begin{array}{rcl}
\mu(X,Y)&=&{1\over{2}}\left\{
\left[{c_0}\left(1-{e^{2\alpha}}\right)
-\de\alpha({{}_\ast}\de\alpha)
\right]\omega(X,Y)+\de\alpha\left(
{\nabla_X}Y+{\nabla_Y}X\right)\right\}\cr
\\
&&-\de\alpha(X)\de\alpha(Y).
\end{array}
\label{12}
\end{equation}

We also have the following PDE deduced from (\ref{3}),
$\forall\,X\in{C^1}(T\cal{M})$
and $\forall\,{\hat f}~\in~{Diff_{loc.}^2}(\cal{M})$:

\begin{equation}
{Tr^1}({\widetilde{\nabla}_X})={Tr^1}({\nabla_X})+n\,\de\alpha(X).
\label{13}
\end{equation}

Thus, we have a serie of PDE deduced from (\ref{1}), in
particular made of  the systems of PDE (\ref{1}), (\ref{10}) and
(\ref{13}). But there are  alternative versions of these PDE in
which the function
$\alpha\in{C^2}(\cal{M},{\R})$ doesn't appear.  These latters are
the following: from the system (\ref{1}), one deduces,
$\forall\,{\hat f}\in{Diff_{loc.}^1}(\cal{M})$:

\begin{equation}
      \left\{
         \begin{array}{ll}
{{\hat f}^\ast}{\hat \omega}={\hat \omega}\\
{\rm and\ with}\,\det(J({\hat f}))\not= 0,
         \end{array}
      \right.
\label{14}
\end{equation}

where ${\hat\omega}=\omega/{|{\det(\omega)}|^{1/n}}$. Then by
prolongation,  with $\widehat{\nabla}$ and $\hat\rho$ being
respectively the connexion of  Levi-Civita and the Riemann
curvature tensor associated to $\hat\omega$,  one obtains
$\forall\,{\hat f}\in{Diff_{loc.}^3}(\cal{M})$:

\begin{equation}
{{}^{\hat{f}}}{\widehat{\nabla}}={\widehat{\nabla}}
\label{15}
\end{equation}

and

\begin{equation}
{{}^{\hat{f}}}{\hat\rho}={\hat\rho}.
\label{16}
\end{equation}

In the latter system (\ref{16}) of PDE (first order), it has to
be noted that
${\hat\rho}=\tau$, i.e. $\hat\rho$ is the Weyl curvature tensor
associated  to the metric $\omega$. Furthermore, with the
assumption that the conformal structure  is flat, one has
${\hat\rho}=0$ and ${\hat\sigma}=0$. But in  general,  it is
noteworthy to add that if $\tau=0$ then obviously $\sigma$ doesn't
vanish. Then, the conformal Lie pseudogroup $\Gamma_{\widehat{G}}$
is   the set of functions ${\hat{f}}\in{Diff_{loc.}^3}(\cal{M})$
satisfying the following  involutive system of PDE:

\begin{equation}
      \left\{
         \begin{array}{ll}
{{\hat f}^\ast}{\hat \omega}={\hat \omega}
\hskip 1cm {\rm and\ with}\hskip 1cm \det(J({\hat f}))\not= 0\\
{{}^{\hat{f}}}{\widehat{\nabla}}={\widehat{\nabla}},
         \end{array}
      \right.
\label{17.1}
\end{equation}

completed with a third system of PDE of order 3 defined
$\forall\,X,Y\in{C^2}(T\cal{M})$ by:

\begin{equation}
{{}^{\hat{f}}}{{\widehat{\nabla}}_X}{{}^{\hat{f}}}{{\widehat{\nabla}}_Y}=
{{\widehat{\nabla}}_X}{{\widehat{\nabla}}_Y}.
\label{17.2}
\end{equation}

This system is formally integrable if and only if $\tau=0$ (from the
Weyl theorem) and involutive because the corresponding symbol
${\widehat{M}}_3$ vanishes. From a terminological point of view, one
shall say that the system  (\ref{17.1})-(\ref{17.2}) is the ``Lie
form" of the system made of the PDE  (\ref{1}), (\ref{10}) and
(\ref{13}) to which one adds the third order  system deduced from
the expression of
${{}^{\hat{f}}}{{\widehat{\nabla}}_X}{{}^{\hat{f}}}{{\widehat{\nabla}}_Y}$
with respect to ${\nabla_X}{\nabla_Y}$ and
$\alpha\in{C^3}(\cal{M},{\R})$.  We shall call this latter system
equivalent to the system (\ref{17.1})-(\ref{17.2}),  the
``deformed" or ``extended" system. It is remarkable that the
deformed  system brings out a supplementary system in comparison
with its Lie form  (\ref{17.1})-(\ref{17.2}). It is about the
system (\ref{10}).  If we consider the second order
sub-system deduced from the system  (\ref{17.1})-(\ref{17.2}), we
notice it is still formally integrable  (again because of the Weyl
theorem) but it is no longer involutive  because the symbol
${\widehat{M}}_2$ of (\ref{17.1})-(\ref{17.2})  is only 2-acyclic.

\addcontentsline{toc}{subsection}{3.2. The conformal  Lie groupo\"{\i}d}
\subsection*{\centering 3.2. The conformal  Lie groupo\"{\i}d}

\hskip 6mm Before defining this groupo\"{\i}d, we need to recall
some definitions concerning the sheafs of the k-jet fiber bundles
(Kumpera {\it et al.} 1972). First of all, we denote ${J_k}(\cal{M})$ the
affine fiber bundle of the k-jets of local
${C^\infty}(\cal{M},\cal{M})$  functions on $\cal{M}$,
$\theta_{\cal{M}}$  (or simply $\theta$) the sheaf of local rings
of  germs of continuous  functions on $\cal{M}$ with values in $\R$,
and $J_k$ the affine  fiber bundle of the k-jets of local functions
in ${C^\infty}(\cal{M},{\R})$. Then in what follows, we  agree to
underline all the names used for the sheafs of germs of  local
continuous sections associated with the various fiber bundles.  Now,
if $\epsilon$ is a sheaf of $\theta$-modules on $\cal{M}$,  we
conventionally define the sheaf $\underline{{J_k}(\epsilon)}$ as:

$$\underline{{J_k}(\epsilon)}\equiv\underline{J_k}
{\otimes_{{}_{\!{\theta}}}}\epsilon.$$

Then, we have the injective sheafs map:

$$\begin{array}{rrcl}
j_k:&\underline{{C^\infty}(\cal{M})}&\longrightarrow &\underline{{J_k}(\cal{M})}\\
\\
&{[m\rightarrow f(m)]_{{}_x}}&\longrightarrow &{[m\rightarrow {j_k}(f)(m)]_{{}_x}}
\equiv{j_k}([f])_{{}_x},
\end{array}$$

where ${j_k}(f)(m)$ is the set of germs at $m$ of the derivatives of $f$
up to order $k$, and $[\,\,\,]_{{}_x}$ obviously being the equivalence
classes of local sections at $x\in \cal{M}$. Let us also denote  
$Aut(\cal{M})$ the sheaf
of germs of functions $f\in{Diff_{loc.}^{\infty}}(\cal{M})$.
The source map such as:

$$\begin{array}{rrcl}
\alpha_k:&{J_k}(\cal{M})&\longrightarrow &\cal{M}\\
&{j_k}(f)(x)&\longrightarrow &x,
\end{array}$$

and the target map:

$$\begin{array}{rrcl}
\beta_k:&{J_k}(\cal{M})&\longrightarrow &\cal{M}\\
&{j_k}(f)(x)&\longrightarrow &f(x),
\end{array}$$

are submersions on $\cal{M}$. One defines the composition on ${J_k}(\cal{M})$ by:

$${J_k}(g)(y)\,.\,{J_k}(f)(x)={J_k}(g\circ{f})(x),$$

with $y=f(x)$. The units in $\underline{{J_k}(\cal{M})}$ are the elements
${j_k}([id]){{}_x}$ and they can  be identified with the points $x\in{\cal{M}}$.
Then, let ${\Pi_k}(\cal{M})$ be the Lie  groupo\"{\i}d of invertible elements of
${J_k}(\cal{M})$. The elements of ${\Pi_k}(\cal{M})$ are the k-jets of the functions
$f\in{Diff_{loc.}^{\infty}}$.
$\underline{{J_k}(\cal{M})}$ (resp.$\underline{{\Pi_k}(\cal{M})}$) is also  
the sheaf of
germs of local continuous sections $f_k$ of $\alpha_k$ (resp.
${\alpha_k}/{\Pi_k}(\cal{M})$). The sheaf map $j_k$ can be also restricted to  
the sheaf map:

$${j_k}:Aut({\cal M})\longrightarrow\underline{{{\Pi}_k}({\cal M})}.$$

An element $[f_k]\in\underline{{\Pi_k}(\cal{M})}$ shall be called
``admissible" if
$f={\beta_k}\circ{f_k}\in{Aut(\cal{M})}$ (i.e.
$\det([{j_1}(f)])\not=0$). The admissible  elements are the germs
of continuous sections $f_k$ of
${\alpha_k}:{\Pi_k}(\cal{M})\rightarrow{\cal{M}}$ such as
${\beta_k}\circ{f_k}\in{Diff_{loc.}^\infty}(\cal{M})$. We denote
${\Gamma_k}(\cal{M})$ the sub-sheaf of admissible elements of
$\underline{{\Pi_k}(\cal{M})}$. Then, we can define the sheaf
epimorphism of groupo\"{\i}ds:

$${j_k}:Aut({\cal M})\longrightarrow{{\Gamma}_k}({\cal M}).$$

Finally, we define the source map $a_k$ and the target map $b_k$ in
$\underline{{J_k}(\cal{M})}$ by:

$$\begin{array}{rrcl}
a_k:&\underline{{J_k}(\cal{M})}&\longrightarrow&\cal{M}\\
&[{\sigma_k}]_{{}_x}&\longrightarrow&x,
\end{array}$$

$$\begin{array}{rrcl}
b_k:&\underline{{J_k}(\cal{M})}&\longrightarrow &\cal{M}\\
&[{\sigma_k}]_{{}_x}&\longrightarrow &{\beta_k}\circ{\sigma_k}(x),
\end{array}$$

and the canonical projection $\Pi{}_q^p$ ($p\geq{q}$) by:

$$\begin{array}{rrcl}
{\Pi{}_q^p}:&\underline{{J_p}(\cal{M})}&\longrightarrow  
&\underline{{J_q}(\cal{M})}\\
&[{f_p}]_{{}_x}&\longrightarrow &[{f_q}]_{{}_x}.
\end{array}$$

Thus, at the sheafs level, the non-linear finite Lie groupo\"{\i}d
$\rtc$ of the conformal pseudogroup $\Gamma_{\widehat{G}}$ is the set
of germs of continuous sections in ${\Gamma_3}(\cal{M})$ satisfying the algebraic
equations over each point $x\in\cal{M}$, obtained by substituting the germs
$[{\hat{f}_3}]_{{}_x}\in{\Gamma_3}(\cal{M})$ for the derivatives of $\hat{f}$ 
up to order three in the system of PDE (\ref{17.1})-(\ref{17.2}). In other  
terms, we
substitute $[{\hat{f}}_3]$ for ${j_3}([\hat{f}])$. More precisely, one
factorizes each differential operator of the system (\ref{17.1})-(\ref{17.2})  
with the
operators $j_k$ $(k=1,2,3)$. Then, one defines the morphisms:

$$
\left\{
\begin{array}{rrcl}
M(\hat{\omega}):&{j_1}({\cal M})&\longrightarrow&{S_2}{T^{\ast}}{\cal M}\\
L({j_1}(\hat{\omega})):&{j_2}({\cal M})&\longrightarrow&
T{\cal M}\otimes\tsm\otimes{J{}_1^{\ast}}(T{\cal M})\\
K({j_2}(\hat{\omega})):&{j_3}({\cal M})&\longrightarrow&
T{\cal M}\otimes\tsm\otimes{J{}_1^{\ast}}(T{\cal M})
\otimes{J{}_2^{\ast}}(T{\cal M})
\end{array}
\right.
$$

by the respective following relations:

$$
\left\{
\begin{array}{rcl}
{\hat{f}^\ast}{\hat{\omega}}&=&M(\hat{\omega})\circ{j_1}(\hat{f})\\
{{}^{\hat{f}}}{\widehat{\nabla}}&=&L({j_1}(\hat{\omega}))
\circ{j_2}(\hat{f})\\
{{}^{\hat{f}}}{\widehat{\nabla}}{{}^{\hat{f}}}{\widehat{\nabla}}&=&
K({j_2}(\hat{\omega}))\circ{j_3}(\hat{f}),
\end{array}
\right.
$$

and consequently $\rtc$ can be rewritten as the system of PDE made
of the two systems of PDE $\rdc$:

\begin{equation}
\left\{
\begin{array}{rcl}
{{}^{[{\hat{f}}_1]}}{\hat{\omega}}\,&{\buildrel{def.}\over{\equiv}}&
\,M({\hat{\omega}})([{\hat{f}}_1])={\hat{\omega}}
{\hskip .2cm}{\rm with}{\hskip .2cm}\det([{\hat{f}}_1])\not=0
{\hskip .2cm}{\rm and}{\hskip .2cm}\det([{j_1}({\hat{f}})])\not=0\cr
{{}^{[{\hat{f}}_2]}}{\widehat{\nabla}}\,&{\buildrel{def.}\over{\equiv}}&
\,L({j_1}({\hat{\omega}}))([{\hat{f}}_2])={\widehat{\nabla}},
\end{array}
\right.
\label{18.1}
\end{equation}

completed with the third system:

\begin{equation}
{{}^{[{\hat{f}}_3]}}({\widehat{\nabla}}{\widehat{\nabla}})\  
{\buildrel{def.}\over{\equiv}}\
K({j_2}({\hat{\omega}}))([{\hat{f}}_3])={\widehat{\nabla}}{\widehat{\nabla}}.
\label{18.2}
\end{equation}

$\Gamma_{\widehat{G}}$ is then the set of germs
$[{\hat{f}}]\in{Aut(\cal{M})}$ such as
${j_3}([{\hat{f}}])\in\rtc$. Now since ${{\widehat{M}}_3}=0$, one
also has the equivalence $\rtc\simeq\rdc$ and in order to  work out
the sophisticated non-linear Spencer complex
of $\rtc$, it is sufficient, as we shall see later, to obtain it
for $\rdc$. That is  because the exactness of this complex of
length two, only needs the symbol  of the corresponding Lie
groupo\"{\i}d to be 2-acyclic. It is precisely  the case for
${\widehat{M}}_2$. Thus, the discussion in what follows will concern
exclusively $\rdc$ defined by the ``deformed" or ``extended"
system:
$\forall\,X,Y\in{C^2}(T\cal{M})$  and
$\forall\,[\alpha_2]=([\alpha],[\beta],[\mu])\in\Jds$,

\begin{equation}
{\rdc}:\left\{
\begin{array}{rcl}
{{}^{[\fuc]}}\omega&=&{e^{2[\alpha]}}\omega
{\hskip .2cm}{\rm with}{\hskip .2cm}\det([\fuc])\not=0
{\hskip .2cm}{\rm and}{\hskip .2cm}\det([{j_1}({\hat{f}})])\not=0\\
{{}^{[\fdc]}}{\nabla_X}Y&=&{\nabla_X}Y+[\beta](X)Y+
[\beta](Y)X-\omega(X,Y)\,{{}_\ast}[\beta]\\
{{}^{[\fuc]}}\sigma(X,Y)&=&\sigma(X,Y)+(n-2)\left([\beta](X)[\beta](Y)-
{{1}\over{2}}\omega(X,Y)[\beta]({{}_\ast}[\beta])\right.\\
&&\left.\mbox{}-[\mu](X,Y)+{{1}\over{2}}[\beta]({\nabla_X}Y+
{\nabla_Y}X) \right).
\end{array}
\right.
\label{19}
\end{equation}

Let us precise again that the equation (\ref{10}) makes sense, and  in
the formula (\ref{19}) the last equation must be considered. Indeed, in
case of a conformally non-flat background  metric $\omega$, the Lie
form of $\widehat{\cal R}_2$ given by the set of equations (\ref{14}),
(\ref{15}) and (\ref{16}), is equivalent to the set of equations made of
the two
first equations in formula (\ref{19}) together with the  equation (\ref{5}).
Then setting $\widehat{\rho}=\tau=0$ (i.e. a conformally flat
metric $\omega$), doesn't change  this equivalence, but in that
case the equation (\ref{5})  becomes equivalent to equation (\ref{10}), and
thus the  ``extended form" of $\widehat{\cal R}_2$  is the formula
(\ref{19}) and we have one equation more than in the Lie form
case. The latter point is rather important to make the difference
between the two forms. Obviously the groupo\"{\i}d ${\cal R}_2\subset\rdc$
of the Poincar\'{e} pseudogroup  corresponds to the case for which
$[{{\alpha}_2}]=0$. The symbol $M_2$ of ${\cal R}_2$ vanishes and
so is involutive, and ${\cal R}_2$ is not formally integrable
unless $\rho_s$  is a constant. The present suggested model is
associated to a particular split exact short sequence of
groupo\"{\i}ds (not of Lie groupo\"{\i}ds because of $\rsd$):

$$1\longrightarrow{{\cal R}_2}\longrightarrow\rdc\longrightarrow\rsd
\longrightarrow{1},$$

and in order to have a relative exact non-linear (even so
fractional!) complex  associated to $\rsd$, the complex associated
to ${\cal R}_2$ would also have to be exact. This is possible only
if ${\cal R}_2$ is formally integrable and consequently involutive
($M_2=0$). Setting these conditions it follows that the relation
(\ref{11}) must be satisfied at the sheaf level. Then from relation
(\ref{11}) and (\ref{19}), one deduces $\rsd$ is the set of
elements $[{\alpha_2}]\in\Jds$ such as $\forall\,{X,Y\in}\,T\cal{M}$
and $\forall\,{c_0}\in{\R}$:

\begin{equation}
{\rsd}:{\hskip 1.5cm}\left\{
\begin{array}{rcl}
[\mu](X,Y)&=&{1\over{2}}\left\{
\left[
{c_0}\left(1-{e^{2[\alpha]}}\right)-
[\beta]({{}_\ast}[\beta])
\right]\omega(X,Y)\right.\\
&&\Bigl.+[\beta]\left(
{\nabla_X}Y+{\nabla_Y}X
\right)
\Bigr\}-[\beta](X)[\beta](Y),
\end{array}
\right.
\label{20}
\end{equation}

and only with these conditions does $\rdc$ reduce indeed to the system chosen
by J.-F.~Pommaret.  Thus, $[\mu]$ is completly defined from
$[\alpha]$ and $[\beta]$  so that the symbol $\msd$ of $\rsd$
obviously vanishes. Hence, $\rsd$ is  involutive and one has the
equivalence:

$$\rsd\simeq\rsu=\Jus,$$

deduced from the embedding of $\Jus$ in $\Jds$ defined by the
system (\ref{20}).  Consequently, one has to work out the complex
associated to $\rsu$ such that:

$$1\longrightarrow{{\cal R}_2}\longrightarrow\rdc\mapright{\phi_0}
\rsu\longrightarrow{1}.$$

The sequences above  are  sequences of group{\"\i}ds but not of  Lie
groupo{\"\i}ds. That $\rsd
\equiv{\underline{J_1}}$ is a  groupo{\"\i}d can be seen directly from the
definition of $\phi_0$ as we shall see in the sequel or first  by considering
locally, above each pair of open subsets $U\times V\subset {\cal M}\times{\cal
M}$ the corresponding associated trivial local groupo{\"\i}ds  $U\times
G_{U\times V}~\times~V\simeq~{\cal R}_{2/U\times V}$ (and with analogous  
expressions for the
other $\cal R$'s). Then we  obtain a corresponding sequence of algebraic
groups on a finite projective free module for the $G$'s. Since they are
algebraic they are splittable and the sequence is split exact. Therefore we
can find a  splitting of groups (such as an Iwasawa decomposition for
instance), by a good choice of back-connection. Then
${\rsu}_{/U\times V}$  can be canonically injected in $\widehat{\cal
R}_{2/U\times V}$ so that it acquires locally an algebraic group structure.
Then, by gluing over all pairs of open subsets
$U\times V$, we deduce the sequence of  groupo{\"\i}ds.
Let us add from the definition of $\rdc$ that $\rsu$ is a natural
bundle associated to $\rdc$. From a physical point
of view, it is important to notice that
$\mu$ may be considered as an Abraham-E\"{o}tvos type tensor,
leading to a first physical interpretation (up to a constant for
units) of $\beta$ as the acceleration 4-vector of gravity. On the
other hand, completly in agreement with J.-F. Pommaret,
$\alpha$ being associated with the dilatations it can be considered
as a relative  Thomson type temperature (again up to a constant for
units):

\begin{equation}
\alpha=\ln({T_0}/T),
\label{21}
\end{equation}

where $T_0$ is a constant temperature of reference associated with
the base  space-time
$\cal{M}$. Then, one can easily define the epimorphism $\phi_0$ by
the relations
$\forall\,[\fdc]\in\rdc$:

\begin{equation}
{\phi_0}([\fdc])=\left\{
\begin{array}{rcl}
[\alpha]&=&{{1}\over{n}}{\ln}|\det({[\fuc]})|\cr
[\beta]&=&{1\over{n}}{Tr_\omega}\left({{}^{[\fdc]}}\nabla-\nabla\right).
\end{array}
\right.
\label{22}
\end{equation}

Thus, we get a first diagram:

\begin{equation}
\begin{array}{ccccc}
&&1&&1\cr
&&\mapdown{}&&\mapdown{}\vespace\cr
1&\mapright{}&{\Gamma_G}&\mapright{j_2}&{{\cal R}_2}\vespace\cr
&&\mapdown{}&&\mapdown{}\vespace\cr
1&\mapright{}&{\Gamma_{\widehat{G}}}&\mapright{j_2}&\rdc\vespace\cr
&&\mapdown{}&&\mapdown{\phi_0}\vespace\cr
1&\mapright{}&\theta&\mapright{j_1}&\Jus\vespace\cr
&&&&\mapdown{}\vespace\cr
&&&&0
\end{array}
\label{23}
\end{equation}

Before presenting the various non-linear Spencer complexes we shall recall
briefly certain definitions and results of this theory. The Spencer Theory
presentation we give below is rather minimal since we  think that it is
impossible to describe it perfectly in few pages. It is especially a matter of
indicating the notations chosen in the  text and we do not claim to make a
full and complete description. Moreover, this theory is presented and taken up
historically with the ``diagonal method" of Grothendieck
(Kumpera {\it et al.} 1972), and the  results and definitions we give don't mention
it. In this method, two copies of the base space $\cal{M}$ are used. The first
one ${\cal{M}}_1$  (the horizontal component) is attributed to the set of
``points" on which the Taylor coefficients of particular Taylor series are
defined, and the second ${\cal{M}}_2$ (the vertical component),  to points on
which these previous series are evaluated. The independence  of this evaluation
with respect to the points chosen in ${\cal{M}}_1$ allow  to deduce the first
Spencer differential operator (linear or non-linear) as an exterior
differential operator on the horizontal component. Then the second Spencer
differential operator and a particular set of  derivations are deduced from the
equivariance of the first Spencer  differential operator with respect to a
particular groupo\"{\i}d action. Thus, the Spencer cohomology can be seen
mainly as an equivariant cohomology on graded sheafs of the diagonal of
${{\cal{M}}_1}\times{{\cal{M}}_2}$. Actually, we give the results and formulas
after the vertical (diagonal isomorphism) projection on the vertical component
of various  diagonal graded sheafs defined on
${{\cal{M}}_1}\times{{\cal{M}}_2}$ and following in parts a formulation given by
J.-F. Pommaret especially when concerning the definition of the differential
bracket. The other definitions are merely the vertical ``translations" of the
definitions given by Kumpera and Spencer applying the so-called ``$\epsilon$"
isomorphism on the diagonal sheafs.


\newtheorem{definition}{Definition}

\addcontentsline{toc}{subsection}{3.3. The first non-linear Spencer  
differential operators}
\subsection*{\centering 3.3.  The first non-linear Spencer differential operators}

\hskip 6mm First, we recall that  $\gtms{}\otith\underline{{J_k}(E)}$,
where $E$ is a
vector bundle, has a natural left $\gtms{}$-module structure.

\begin{definition}
\begin{description}
\item{a)} The linear Spencer operator $D$ is the unique differential
operator ($\R$-linear sheaf map; $k,s\geq{0}$)\par
$$D:\gtms{}\otith\underline{{J_{k+1}}(E)}\longrightarrow\gtms{}\otith
\underline{{J_k}(E)},$$
satisfying the three following conditions:\par
\begin{enumerate}
\item $D\circ{j_{k+1}}=0$,
\item $\forall{\tau_{k+1}}\in\gtms{}\otith\underline{{J_{k+1}}(E)}$,
{\hskip 5mm}$D(\omega\wedge{\tau_{k+1}})=\de\omega\wedge{\tau_k}+
(-)^{{d^0}\omega}\omega\wedge{D{\tau_k}}$ \par where $\omega\in\gtms{}$
is any homogeneous differential form and $\de$ being the exterior differential,
\item $D$ restricted to ${J_{k+1}}(E)$ satisfies:
$${\epsilon_1}\circ{D}={j_1}\circ{\Pi{}_k^{k+1}}-{id_{{J_{k+1}}(E)}},$$
\par
where $\epsilon_1$ is the monomorphism defined by the short
exact sequence:
$$0\longrightarrow{{T^\ast}\!\!{\cal M}\otimes{{J_k}(E)}}
\mapright{\epsilon_1}{{J_1}({{J_k}(E)})}\mapright{\Pi{}_0^1}{{J_k}(E)}
\longrightarrow{0}.$$
\end{enumerate}
\item{b)} The restriction of ($-D$) to the symbol $\gtms{}\otith
\underline{{S_k}{T^\ast}\!\!{\cal M}}\otith\underline{E}$ defines
the $\gtms{}$-linear Spencer map $\underline{\delta}$ such that
$$\underline{\delta}(\omega\wedge{\tau_k})=(-)^{{d^0}\omega}
\omega\wedge\underline{\delta}({\tau_k}),$$\par
with $\omega$ and $\tau_k$ as in a)-3..\car
\end{description}
\label{def1}
\end{definition}

\begin{definition}
We define the r-th Spencer sheaf ($r\geq{1}$) of ${J_k}(E)$, the quotient sheaf
$${{\cal C}{}_k^r}(E)=\gtms{r}\otith\underline{{J_k}(E)}/{\zeta_k}\circ
\underline{\delta}(\gtms{r-1}\otith\stms{k+1}\otith\underline{E}),$$
where $\zeta_k$ is the monomorphism defined by the short exact sequence:
$$0\longrightarrow\stm{k}\otimes{E}\mapright{\zeta_k}{J_k}(E)
\mapright{\Pi{}_{k-1}^k}{J_{k-1}}(E)\longrightarrow{0}.$$\car
\label{def2}
\end{definition}

Let us add that ${{\cal C}{}_k^r}(E)$ has a module structure on the  
$\theta$-algebra
$\underline{J_k}$ and we set ${{\cal C}{}_k^r}(E)=0$ for $r>n$ and
${{\cal C}{}_k^0}(E)=\underline{{J_k}(E)}$.

\begin{definition}
 The operator $D$ can be factorized with a right $\theta$-linear operator $D'$ on
${{\cal C}{}_k^r}(E)$:
$$D':{{\cal C}{}_k^r}(E)\longrightarrow{{\cal C}{}_k^{r+1}}(E).$$\car
\label{def3}
\end{definition}

Contrary to the operator $D$, there is no loss of order on $k$. More precisely, 
$D'$ is such that the following diagram of split exact sequences is commutative:

$$\begin{array}{ccc}
0&&0\cr
\mapdown{}&&\mapdown{}\vespace\cr
\stms{k+1}\otith\underline{E}&\mapright{(-\delta)}&
\underline{\tsm}\otith\stms{k}\otith\underline{E}\vespace\cr
\mapdown{\zeta_{k+1}}&&\mapdown{id\otimes{{\zeta}_k}}\vespace\cr
\underline{{J_{k+1}}(E)}&\mapright{D}&\tsms\otith\underline{{J_k}(E)}\vespace\cr
\mapdown{\Pi{}_k^{k+1}}&&\mapdown{\rho{}_k^{k+1}}\vespace\cr
\underline{{J_k}(E)}&\mapright{D'}&{{\cal C}{}_k^1}(E)\vespace\cr
\mapdown{}&&\mapdown{}\vespace\cr
0&&0\vespace
\end{array}$$

It is to be noted that the sequences being split, then $D'$ is built up from a 
connexion

$$c{}^k_{k+1}:\underline{{J_k}(E)}\longrightarrow\underline{{J_{k+1}}(E)},$$

such that by definition

$${\Pi{}^{k+1}_k}\circ{c{}^k_{k+1}}={id_k}.$$

But quotienting, then by definition, $D'$ is independent of the choice of connexion
$c{}^k_{k+1}$. Hence, whatever is $c{}^k_{k+1}$, one has
$$D'={\rho{}_k^{k+1}}\circ{D}\circ{c{}^k_{k+1}}.$$
Finally, for $r\geq{1}$, these definitions can be extended to the tangent  
bundle $R_k$ of
${\cal R}_k$ instead of ${J_k}(E)$, and to its corresponding symbol $M_k$. But
in this case, from the definition of $D'$, 1) ${\cal R}_{k+1}$ must be a  
fibered manifold,
2) to make a choice of connexion  $c{}^k_{k+1}$ we must have the epimorphism
${{\cal R}_{k+1}}\longrightarrow{{\cal R}_k}\longrightarrow{0}$, and 3) the system 
${\cal R}_k$ must be formally transitive, i.e. we must have the epimorphism
${{\cal R}_k}\longrightarrow{\cal M}\longrightarrow{0}$. Also we shall use  
the definitions
and notations:

$${{\cal C}{}_k^r}=\gtms{r}\otith{R_k}/{\zeta_k}\circ
\underline{\delta}(\gtms{r-1}\otith{M_{k+1}}).$$

\begin{definition}
Let us define ${B{}_k^r}({\cal M})$ and ${B{}_k}({\cal M})$ such that
$${B{}_k^r}({\cal M})=\gtms{r}\otith\underline{J_k},$$ and
$${B{}_k}({\cal M})={\buildrel{r}\over\oplus}{B{}_k^r}({\cal M}).$$\car
\label{def4}
\end{definition}

${B{}_k}({\cal M})$ has a natural structure of  graded  
$\underline{J_k}$-algebra  defined
by the exterior product of the $J_k$-valued forms, since we have the equivalence:

$${B{}_k}({\cal M})\equiv\gtms{}\otith\underline{J_k}.$$

${B{}_k}({\cal M})$ inherits also a natural structure of left graded  
$\gtm{}$-module  where
the external operation on ${\buildrel{r>0}\over\oplus}{B{}_k^r}({\cal M})$ is  
 the exterior
product $\omega\wedge\mu$ of $\omega\in\gtm{}$  and $\mu\in{B{}_k}({\cal M})$  
with respect
to the pairing on ${B{}_k^0}({\cal M})$:

$$([f],[{g_k}])\in\theta\times\underline{J_k}\longrightarrow
{j_k}([f]).[g_k]\in\underline{J_k}.$$

Let us denote by ${Der_a}{B{}_k}({\cal M})$ the sheaf of germs of  
``admissible" graded
derivations $\cal D$ of ${B{}_k}({\cal M})$, i.e. if $\cal D$ is of degree  
$p$, one has:

$$
{\cal D}(\underline{J{}_k^0})\subset\gtms{p}\otith\underline{J{}_k^0},
$$

$$
{\cal D}(\gtms{r})\subset\gtms{r+p},
$$

where $\underline{J{}_k^0}$ is the kernel of the target map $\beta_k$ defined on 
$\underline{J{}_k}$:

$${\beta_k}:\underline{J{}_k}\longrightarrow\theta.$$

Obviously ${Der_a}{B{}_k}({\cal M})$ is endowed with the bracket of  
derivations $[\quad,\quad]$,
i.e. if  ${\cal D}_i$ is a derivation of degree $p_i$, then  we have:

$$
[{\cal D}_1,{\cal D}_2]={\cal D}_1\circ{\cal D}_2 -
(-)^{p_1 p_2}{\cal D}_2\circ{{\cal D}_1}.
$$

Finally, we also have

$$
{Der_a}{B{}_k}({\cal M})={\buildrel{r}\over\oplus}{Der{}_a^r}{B{}_k}({\cal M}),
$$

where ${Der{}_a^r}{B{}_k}({\cal M})$ is the module of admissible derivations  
of degree
$r$ on ${B{}_k}({\cal M})$.

\begin{definition}
One defines $\twinu$ (the ``twisting" of $\de$), the non-linear differential
operator
$$\twinu:{\Gamma_{k+1}}{\cal M}\longrightarrow{Der{}_a^1}{B{}_k}({\cal M})$$
such that:
$$\twinu[f_{k+1}]=\de - Ad[f_{k+1}]\circ{\de}\circ{Ad[f{}^{-1}_{k+1}]},$$
where $Ad[f_{k+1}]$ is the contravariant action at the sheafs level of  
$[f_{k+1}]$ on
$\gtms{}\otith\underline{J_k}$ corresponding  to the action of the pull-back of 
$f\in{Aut({\cal{M}})}$ on the tensors of $\gtm{}\otimes{J_k}$.\car
\label{def5}
\end{definition}

To this action on $\gtm{}\otimes{J_k}$ corresponds simultaneously
an action of  the ``pull-back-push-forward of
$f\in{Aut({\cal{M}})}$ " on the tensors of
$\gtm{}\otimes{J_k}({\cal{M}})$. Also, we deduce and define at the
sheaf level,  the extension of $Ad[f_{k+1}]$ on
$\gtms{}\otith\underline{{J_k}({\cal{M}})}$.

\begin{definition}
Let
$$\twinup:{\Gamma_k}{\cal M}\longrightarrow{{\cal C}{}_k^1}(T{\cal{M}}),$$
 be ``the first non-linear Spencer operator" such that
$\forall[f_k]\in{\Gamma_k}{\cal{M}}$:
\begin{equation}
\twinup([f_k])={\rho{}_k^{k+1}}\circ\left[
\de - Ad[f_{k+1}]\circ{\de}\circ{Ad[f{}^{-1}_{k+1}]}\right](id_k),
\label{24}
\end{equation}
where $[f_{k+1}]={c{}_{k+1}^k}([f_k])$ and $id_k\in\underline{{J_k}({\cal{M}})}$ is
the prolongation up to order $k$ of $id_{\cal{M}}~\equiv~{id_0}$.\car
\label{def6}
\end{definition}

Sometimes this definition is given using the functor $j_1$ instead of $\de$.
The result is that the derivation is made with respect to the source and not
the target and difficulties appear when defining the brackets given further.
Then Spencer defined the isomorphism $ad$ of degree zero:

\begin{definition}
One defines the isomorphism $ad$, the operator
$$
ad:{{{\cal C}{}_{k+1}^{\bullet}}(T{\cal{M}})}
\longrightarrow{Der_a}{B_k}({\cal M})$$
such that $\forall{v}\in{B_k}(\cal{M})$ and
$\forall{u_{k+1}}\in{{{\cal C}{}_{k+1}^{r}}(T{\cal{M}})}$:
\begin{enumerate}
\item $ad({{\cal C}{}_{k+1}^{\bullet}}(T{\cal M})){\buildrel{def.}\over{\equiv}}
{Der_\Sigma}{B_k}({\cal  
M})={Der_{\Sigma,k}^{\bullet}}\subseteq{Der_a}{B_k}({\cal M})$,
\item ${\cal L}:\gtms{}\otith\underline{{J_k}(T{\cal M})}\longrightarrow
{Der_a}{B_k}({\cal M})$ being the Lie derivative of degree zero such that
$${\cal  
L}(u_k)v=[i({u_k}),\de]v\equiv{u_k}\wedgebar{\de{v}}+{(-)^r}\de({u_k}\wedgebar{v}),$$
where $\wedgebar$ is the extended Fr\"{o}licher-Nijenhuis product and
$i$ the interior product:
$$i({u_k})v\equiv{u_k}\wedgebar{v},$$
\item $ad(u_{k+1})v=\left[ {\cal{L}}(u_k)+
(-)^{r+1}D'(u_{k+1})\wedgebar\right]v$.
\end{enumerate}\car
\label{def7}
\end{definition}

Initially, $ad$ is the $\theta$-linear sheafs map corresponding to the  
($k+1$)-st order
differential operator ${\cal L}\circ{j_k}$ where ${\cal  
L}:{J_k}(T)\longrightarrow{Der_a}
{B_k}({\cal M})$ is the Lie derivative of degree zero and of order 1. Then,  
one obtains the
first step of the sophisticated non-linear Spencer complexes associated to  
the resolutions
${{\cal C}{}_{k+1}^{\bullet}}(T{\cal M})$ or ${Der_{\Sigma,k}^{\bullet}}$:


$$\begin{array}{ccccccc}
1&\longrightarrow&Aut({\cal M})&\mapright{j_{k+1}}&{\Gamma_{k+1}}{\cal M}&
\mapright{\twinup}&{{{\cal C}{}_{k+1}^1}(T{\cal{M}})}\cr
&&{\Big|}{\Big|}&&{\Big|}{\Big|}&&\mapdown{ad}\vespace\cr
1&\longrightarrow&Aut({\cal M})&\mapright{j_{k+1}}&{\Gamma_{k+1}}{\cal M}&
\mapright{\twinu}&{Der_{\Sigma,k}^1}\vespace
\end{array}$$


where the two horizontal sequences are split exact at
${\Gamma_{k+1}}{\cal M}$. From this, we can determine the induced
fractional differential operator $\twinsu$  such that the following
diagram is commutative (${M_3}=\mc{3}=0$,
$\msd=0$):

\begin{equation}
\begin{array}{ccccccccl}
&&1&&1&&0&&\cr
&&\mapdown{}&&\mapdown{}&&\mapdown{}&&\vespace\cr
1&\longrightarrow&\Gamma_G&\mapright{j_2}&{\cal R}_2&\mapright{\twinup}&
{\cal C}{}_2^1&=&\tsms\otith{R_2}\vespace\cr
&&\mapdown{}&&\mapdown{}&&\mapdown{}&&\vespace\cr
1&\longrightarrow&\Gamma_{\widehat G}&\mapright{j_2}&{\widehat{\cal R}}_2&
\mapright{\twinupc}&\widehat{{\cal C}}{}_2^1&=&\tsms\otith{\widehat{R}_2}\vespace\cr
&&\mapdown{}&&\mapdown{\phi_0}&&\mapdown{\phi_1}&&\vespace\cr
1&\longrightarrow&\theta_{\cal M}&\mapright{j_1}&\underline{J_1}&
\mapright{\twinsu}&\not{\!{\cal C}}{}_1^1&=&B{}_1^1({\cal M})\vespace\cr
&&&&\mapdown{}&&\mapdown{}&&\vespace\cr
&&&&1&&0&&\vespace
\end{array}\label{24bis}
\end{equation}


Obviously, one will determine also $\phi_1$ and in the sequel, one
will call
$B{}_1^1({\cal{M}})$ the sheaf of electromagnetic and gravitational
gauge  potentials.

\addcontentsline{toc}{section}{4. The potentials of interaction and the metric}
\section*{\centering 4. The potentials of interaction and the metric}

\addcontentsline{toc}{subsection}{4.1. The electromagnetic and gravitational  
potentials}
\subsection*{\centering 4.1. The electromagnetic and gravitational potentials}

\hskip 6mm In order to lighten the presentation of the results,
first let us consider the following notations:

\begin{enumerate}
\item[1)] one will merely write $f_k$ and ${j_k}(f)$ instead of respectively  
$[f_k]$ and
${j_k}([f])$,
\item[2)] one will denote  ${T^q}f$ the restriction in $\stms{q}\otith
\underline{T{\cal{M}}}$ of $[f_p]\in\underline{{J_p({\cal{M}})}}$
($p\geq{q}\geq{1}$),
\item[3)] one will denote  $T({T^r}f)$ the restriction in
$\tsms\otith\stms{r}\otith\underline{T{\cal{M}}}$ of
${j_1}([f_s])\in\underline{{J_1}({J_s({\cal{M}})})}$ ($s\geq{r}\geq{0}$),
\item[4)] one sets analogous conventions concerning the presence of ordinary
parenthesis in the notations for the differential of germs of
the tangent maps
$d{T^q}f$ and $d({T^q}f)$ (let us note in these notations that $d$ is not the 
exterior differential that one denotes by $\de$, but stands for
differential maps).
\end{enumerate}

Then, one has the following set of results $\forall\,\fdc\in\rdc$ and
$\forall\,{X}\in{\cal C}^1(T{\cal{M}})$:

\begin{eqnarray}
{Tr^1}({{}^{{j_1}(\fuc)}}{\nabla_X})&=&{{1}\over{2}}{\sum^n_{i,j=1}}
{\tilde{\omega}^{i,j}}\left\{<<d(\omega){\circ}\hat{f}{|}T(\hat{f}).X>
{|}T\hat{f}.{e_i}{\otimes}T\hat{f}.{e_j}>\right.\nonumber\\
&&\left.\mbox{}+2<\omega{\circ}\hat{f}{|}<d(T\hat{f})
{|}X>.{e_i}{\otimes}T\hat{f}.{e_j}>\right\}\label{25.1}\\
&=&{Tr^1}({\nabla_X})+n<\de(\alpha){|}X>,
\label{25.2}
\end{eqnarray}

and

\begin{eqnarray}
{Tr^1}({{}^{\fdc}}{\nabla_X})&=&{{1}\over{2}}{\sum^n_{i,j=1}}
{\tilde{\omega}^{i,j}}\left\{<<d(\omega){\circ}\hat{f}{|}T\hat{f}.X>
{|}T\hat{f}.{e_i}{\otimes}T\hat{f}.{e_j}>\right.\nonumber\\
&&\left.\mbox{}+2<\omega{\circ}\hat{f}{|}<dT\hat{f}
{|}X>.{e_i}{\otimes}T\hat{f}.{e_j}>\right\}\label{26.1}\\
&=&{Tr^1}({\nabla_X})+n<\beta{|}X>.\label{26.2}
\end{eqnarray}

Now, let $\widehat{\chi}^{(2)}$ be an element of  
$\underline{{J_1}({J_2(T{\cal{M}})})}$,
and its components $\widehat{\chi}_q$ ($q=0,1,2$), i.e. the restrictions of
$\cc{2}$ to $\tsms\otith\stms{q}\otith\underline{T{\cal{M}}}$ such that
$\forall\ftc\in\rtc$:

\begin{equation}
\cc{2}=\ftci{3}{\circ}{j_1}(\fdc)-id_3,\label{27}
\end{equation}

where by abuse of notations $id_3$ is the image of $id_3\in{J_3}(T{\cal M})$ by 
the canonical injection $\displaystyle{J_3}(T{\cal M})\longrightarrow
{J_1}({J_2}(T{\cal M}))$. In particularly, one has the relation
(Pommaret 1994):

$$\begin{array}{rcl}
\ccc{0}&=&\ftci{1}\circ{j_1}(\hat{f})-id_1\\
&\equiv&{\widehat{A}}-id_1.
\end{array}$$

It follows that $\twinupc$ satisfies:

\begin{equation}
{\epsilon_1}\circ\twinupc(\fdc)\equiv\tc{2}=\cc{2}\circ(\widehat{B}\otimes{id_2}),\label{28}
\end{equation}

with $\widehat{B}={\widehat{A}^{-1}}$ and $\ftc={c{}_3^2}(\fdc)\in\rtc$.
In particular, $\tcc{0}$ and $\tcc{1}$ satisfy the relations:

$$\tcc{0}=id_{T{\cal{M}}}-T\hat{f}\circ{T}(\hat{f})^{-1}\equiv{id_{T{\cal{M}}}}
-\widehat{B}=\ccc{0}\,.\,\widehat{B}\in{T{\cal{M}}}\otimes\tsm,$$

and ${\forall}\,X,Y\in{T{\cal{M}}}$ (Pommaret 1994):

$$
<dT\hat{f}|{X}>\circ\,\tcc{0}\,.\,Y+T\hat{f}\circ<\tcc{1}{|}Y>.\,X=
<d(T\hat{f})-dT\hat{f}|\widehat{B}.Y>.\,X,
$$

where $<dT\hat{f}|{X}>$, $<\tcc{1}{|}Y>\equiv<\ccc{1}{|}\widehat{B}.Y>$
and $<d(T\hat{f})-dT\hat{f}|\widehat{B}.Y>$,
are considered as elements of ${T{\cal{M}}}\otimes\tsm$.
Then, if one substitutes

$$
{T\hat{f}^{-1}}\circ{T(\hat{f})}.X=(\ccc{0}+id).X
$$

for $X$ in the relations (\ref{26.1})-(\ref{26.2}), and also considering the  
Schwarz equalities:

$$
<dT\hat{f}|{X}>.\,Y=<dT\hat{f}|{Y}>.\,X,
$$

one obtains by subtracting the result from the relations (\ref{25.1})-(\ref{25.2}) 
$\forall\,X\in{T{\cal{M}}}$:

\begin{equation}
{{1}\over{n}}Tr^1\left[<\ccc{1}{|}X>+{\nabla_{\ccc{0}.X}}\right]=
<\de(\alpha)-\beta|{X}>-<\beta|{\ccc{0}.X}>.\label{29}
\end{equation}

From this latter relation, one can define the electromagnetic potential vector 
$\pve\in\tsms$ by $\forall\,{X}\in\underline{T{\cal{M}}}$:

\begin{equation}
<\pve|{X}>={{1}\over{n}}Tr^1\left[<\tcc{1}{|}X>+
{\nabla_{\tcc{0}.X}}\right]=
<\de(\alpha)|{\widehat{B}.X}>-
<\beta|{X}>.\label{30}
\end{equation}

In an orthonormal system of coordinates, the latter definition becomes
($i,j,k=1,...,n$):

$$
{\pve_i}={{1}\over{n}}\left({\widehat{\tau}{}^k_{k,i}}+
{\widehat{\tau}{}^k_{\,,i}}{\gamma{}^j_{j\,k}}
\right)={\widehat{B}{}^k_i}{\partial_k}\alpha-{\beta_k},
$$

where $\gamma$ is the Riemann-Christoffel 1-form associated to $\omega$ and  
thus satisfying:

$$
Tr^1(\nabla_X)=Tr^1(\gamma(X)).
$$

Prolonging the relation (\ref{29}) (one does not prolong the
relation  (\ref{30}) because $\tcc{2}$ is not the first
prolongation of $\tcc{1}$, contrary to $\ccc{2}$), one
deduces and defines the mixed tensor potential  of gravitation and
electromagnetic $\ptge\in{\stackrel{2}{\otimes}}\tsms$ such  as
$\forall\,{X,Y}\in\tsm$:

\begin{eqnarray}
<\ptge|{Y}\otimes{X}>&=&
{{1}\over{n}}Tr^1\left[
{i_Y}<\tcc{2}|{X}>+<d(\gamma)|{Y}\otimes\tcc{0}.X>
\right.\nonumber\\
&&\left.\mbox{}+\gamma(<\tcc{1}|{X}>.\,Y)\right]\nonumber\\
&=&<\de(\beta)|\widehat{B}.X\otimes{Y}>-
<\beta|<\tcc{1}|{X}>.\,Y>\nonumber\\
&&\mbox{}-<\mu|{X}\otimes{Y}>\label{31}
\end{eqnarray}

where $i_Y$ is the interior product by $Y$ and $\tcc{2}$ satisfies
$\forall\,\ftc\in\rtc$ and $\forall\,{X,Y,Z}\in{T{\cal{M}}}$:

$$
\begin{array}{l}
<<d(dT\hat{f})-{d^2}T\hat{f}|\widehat{B}.X>|{Y}\otimes{Z}>=\\
\\
{\hskip 3cm}<{d^2}T\hat{f}|{Y}\otimes{Z}>\circ\,\tcc{0}\,.\,X+
<dT\hat{f}|{Y}>\circ<\tcc{1}|{X}>.\,Z\\
{\hskip 3cm}\mbox{}+<dT\hat{f}|{Z}>\circ<\tcc{1}|{X}>.\,Y+
T\hat{f}\circ<<\tcc{2}|{X}>|{Y}\otimes{Z}>.
\end{array}
$$

Again, in an orthonormal system of coordinates ($i,j,k,h=1,...,n$):

$$
{\ptge_{j,i}}={{1}\over{n}}\left({\widehat{\tau}{}^k_{k\,j,i}}+
{\widehat{\tau}{}^k_{j\,,i}}{\gamma{}^h_{h\,k}}+
{\widehat{\tau}{}^k_{\,,i}}({\partial_k}{\gamma{}^h_{h\,j}})
\right)={\widehat{B}{}^k_i}{\partial_k}{\beta_j}-{\mu_{i\,j}}-
{\widehat{\tau}{}^k_{j\,,i}}{\beta_k}.
$$

This definition for $\ptge$ hasn't been determined by J.-F.
Pommaret (1994), he also gave  a different definition for
$\pve$ rather associated to the relation (\ref{29}). In conclusion
to this chapter, one notices  that $\twinsu$ (remaining
to explicit) appears to be a Fr\"{o}benius-type operator and
depends on $\tc{1}$ itself depending on $\pve$ and $\ptge$ as we
shall see further. Lastly, $\phi_1$ is quite defined by the
relations (\ref{30}) and (\ref{31}),  and  if $n=4$, one has 20
scalar gauge potentials. On the other hand, one  can see as well
that the definitions of $\pve$ and $\ptge$ can be deduced from  the
conformal Killing equations on $\underline{{J_1}(T{\cal{M}})}$,
namely
$\forall\,{\xi_{(1)}}\equiv(\xi_0,\xi_1)\in\underline{{J_1}(T{\cal{M}})}$
and $\forall\,\eta\in\theta$ then one has:

$$
{K_0}({\xi^{(1)}})\stackrel{def.}{\equiv}{{1}\over{n}}{Tr^1}(\xi_1+
\gamma(\xi_0))=\eta,
$$

where $K_0$ is the conformal Killing operator. If $K_1$ is its first
prolongation then seting ${K^{(1)}}\equiv(K_0,K_1)$ one obtains obviously from 
the equation above $\forall\,{X}\in\underline{T{\cal{M}}}$ and
$\forall\,\tc{2}\in~{\widehat{\cal{C}}{}^1_2}({T{\cal{M}}})$:

$$
{\phi_1}(\tc{2})(X)={K^{(1)}}(<\tc{2}|{X}>).
$$

\addcontentsline{toc}{subsection}{4.2. The morphisms $\phi_1$, $\twinsu$ and  
the metric of the ``gauge space-time"}
\subsection*{\centering 4.2. The morphisms $\phi_1$, $\twinsu$ and the metric  
of the ``gauge space-time"}

\hskip 6mm From the preceding chapter, obviously one easily notices
by definition  that ${\cal{C}}{}^1_2$ is the kernel of $\phi_1$. On
the contrary $\twinsu$,  depending on $\tc{1}$, is rather more
tricky to determine. The sequences being  split and $\phi_1$ being
$\theta$-linear, one can deduce the important relation that
$\tc{2}$ can be rewritten in the following form:

\begin{equation}
\tc{2}=\tpc{2}+<\pve|\Omega^{(2)}>+<\ptge|{\cal K}^{(2)}>,
\label{31bis}
\end{equation}

where $\tpc{2}\in{\cal{C}}{}^1_2$ and $\Omega^{(2)}$ and ${\cal
K}^{(2)}$ are elements of
$\gtms{}\otith({\widehat{\cal{C}}{}^1_2}\Big/{\cal{C}}{}^1_2)$,
one calls the linear susceptibilities of the vacuum or the
space-time  associated respectively to $\pve$ and $\ptge$. Then, as
a definition of
$\twinsu$, one obtains the following kind of relations:

$$\begin{array}{rcccl}
<a_1|\pve>+<b_1|\ptge>&=&D'\alpha&\equiv&\de\alpha-\beta,\\
<a_2|\pve>+<b_2|\ptge>&=&D'\beta&\equiv&\de\beta-\mu,
\end{array}$$

where $a_1$ is affine with respect to $\de\alpha$, $b_1$  linear
with respect  to $\de\alpha$, $a_2$  linear with respect to
$\de\beta$ and $\beta$, and $b_2$ affine with respect to $\de\beta$
and $\beta$. Thus, $\twinsu$ is fractional  with respect to
$\de\alpha$, $\de\beta$, $\alpha$ and $\beta$ and finally  only
depends on the susceptibilities of $\tc{2}$. Then, the exactness
condition at $\Jus$:

$$
\twinsu\circ{j_1}=0,
$$
the relations (\ref{30})-(\ref{31}) and the commutativity of the diagram
(\ref{24bis}) involve the necessary condition which must be satified by
$\tpc{2}$:

$$\tau_0=0.$$

There remains the metric $\nu$ of the ``gauge (or observable or measurable)
space-time" defined in considering $\widehat{B}$ as a field of tetrads
$\forall\,{X,Y}\in{T{\cal{M}}}$:

$$<\nu|{X}\otimes{Y}>=<\omega|\widehat{B}.X\otimes\widehat{B}.Y>.$$

Thus, one has the general relation between $\nu$ and $\omega$:

$$\nu=\omega+\mbox{linear\,and\,quadratic\,terms\,in\,$\pve$
\,and\,$\ptge$}.$$

Then from this metric $\nu$, one can deduce the Riemann and Weyl
curvature  tensors of the gauge space-time. One has  a non-metrical
theory for  the gravitation in the gauge space-time, since clearly
$\nu$  doesn't appear  as a gravitational potential. The space-time
terminology we use is quite  natural in the sense that one has
simultaneously two types of space-time. The  first one, which we call the
``underlying"  or ``substrat" space-time, is endowed  with the metric $\omega$ and
is of constant scalar curvature. The other one,  called the ``gauge
(observable or measurable) space-time", endowed with  the metric
$\nu$, is defined for any scalar curvature and by the gauge
potentials $\pve$ and $\ptge$. It can be considered as the
underlying space-time  deformed by the gauge potentials and the
Weyl curvature does not necessarily vanish, contrary to Pommaret's
assertions (Pommaret 1994, see page 456 and Pommaret 1989). Moreover, from
a continuum mechanics of  deformable bodies point of view, the
metric $\nu$ can be interpreted as the tensor of deformation of the
underlying space-time (Katanaev {\it et al.} 1992, Kleinert 1989).

\addcontentsline{toc}{section}{5. The fields of interaction}
\section*{\centering 5. The fields of interaction}

\addcontentsline{toc}{subsection}{5.1. The second non-linear Spencer operator}
\subsection*{\centering 5.1. The second non-linear Spencer operator}

\hskip 6mm Before giving an explicit expression for this operator
in the  complex of $\Jus$, one will briefly recall its definition,
but before that one needs  a few other definitions
(Pommaret 1989, Kumpera {\it et al.} 1972).

\begin{definition}
Let the ``algebraic bracket" $\{\quad,\quad\}$ on
${J_{k+1}}(T{\cal M})$, be the $\R$-bilinear map ($\forall{k}\geq{0}$)
$$
\{\quad,\quad\}:{J_{k+1}}(T{\cal M}){\times_{{}_{\cal M}}}{J_{k+1}}(T{\cal M})
\longrightarrow{J_k}(T{\cal M})
$$
such as $\forall{\xi_{k+1}},{\eta_{k+1}}\in{J_{k+1}}(T{\cal M})$
one has:
$$
\{{\xi_{k+1}},{\eta_{k+1}}\}=\{{\xi_1},{\eta_1}\}_{{}_k},
$$
where $\{{\xi_1},{\eta_1}\}\in{T{\cal M}}$ is the usual Lie bracket
defined on ${J_1}(T{\cal M})$ and $\{{\xi_1},{\eta_1}\}_{{}_k}$
its lift in ${J_k}(T{\cal M})$.\car
\end{definition}

\begin{definition}
One calls ``differential Lie bracket" on $\underline{{J_k}(T{\cal M})}$, the
bracket $\grobral\quad,\quad\grobrar$ such that:
$$
\grobral\quad,\quad\grobrar:
\underline{{J_k}(T{\cal M})}{\times_{{}_{\cal M}}}
\underline{{J_k}(T{\cal M})}\longrightarrow\underline{{J_k}(T{\cal M})},
$$
and $\forall{\xi_k},{\eta_k}\in\underline{{J_k}(T{\cal M})}$ then
$$
\grobral{\xi_k},{\eta_k}\grobrar=
\{{\xi_1},{\eta_1}\}+
{i_{\xi_0}}D{\eta_{k+1}}-
{i_{\eta_0}}D{\xi_{k+1}},
$$
where $i$ is the usual interior product and $\xi_{k+1}$ and $\eta_{k+1}$
are any lifts of $\xi_k$ and $\eta_k$ in
$\underline{{J_{k+1}}(T{\cal M})}$.\car
\end{definition}

\begin{definition}
For any decomposable elements
$$\alpha=u\otimes{\xi_k}\in\gtms{r}\oplus\jtm{k},$$
$$\beta=v\otimes{\eta_k}\in\gtms{s}\oplus\jtm{k},$$
and defining on $\gtms{}\otith\jtm{k}$ the interior product $i$
by the relation:
$\forall{w}\in\gtms{}$,
$$
{i_\alpha}w=u\wedge{i_{\xi_k}}w,
$$
then with $\de$ being the exterior derivative and
\begin{enumerate}
\item[1)] ${\cal L}$ the Lie derivative on $\gtms{}\otith\jtm{k}$
such that:
$${\cal L}(\alpha)={i_\alpha}\circ\de+
(-)^r\de\circ{i_\alpha},$$
\item[2)] $ad(\alpha)={\cal L}(\alpha)+(-)^{r+1}{i_{{}_{D\alpha}}},$
\end{enumerate}
one defines the ``twisted" bracket on $\gtms{}\otith\jtm{k}$
($\theta$-bilinear),\quad $\grobral\quad,\quad\grobrar$ by the relation:
$$
\grobral\alpha,\beta\grobrar=
[ad(\alpha)v]\otimes{\eta_k}-(-)^{rs}
[ad(\beta)u]\otimes{\xi_k}+
(u\wedge{v})\otimes\grobral{\xi_k},{\eta_k}\grobrar
\in\gtms{r+s}\otith\jtm{k}.
$$\car\label{def10}
\end{definition}

This bracket defines a graded Lie algebra structure on
$\gtms{}\otith\jtm{k}$.

\begin{definition}
One calls $2^{nd}$ non-linear Spencer operator of the
resolution $Der{}^{\bullet}_{\Sigma,k}$, the differential
operator $\twind$ such that:
$$
\twind:{Der{}^1_{\Sigma,k}}\longrightarrow{Der{}^2_{\Sigma,k}},
$$
and $\forall{u}\in{Der{}^1_{\Sigma,k}}$
$$
\twind{u}=[\de,u]-{1\over{2}}[u,u],
$$\car
\end{definition}

\begin{definition}
To this operator $\twind$ corresponds the $2^{nd}$ non-linear
Spencer operator $\twindp$ of the resolution
${{\cal C}^{\bullet}_k}(T{\cal M})$ such that:
$$
\twindp:{{\cal C}^1_k}(T{\cal M})\longrightarrow
{{\cal C}^2_k}(T{\cal M}),
$$
and $\forall\tau\in{{\cal C}^1_k}(T{\cal M})$
$$
\twindp\tau=D'\tau-{1\over{2}}\grobral\tau,\tau\grobrar',
$$
where $\grobral\quad,\quad\grobrar'$ is the quotient twisted
bracket.\car
\end{definition}

\begin{definition}
With the latter definitions, the ``sophisticated non-linear
Spencer complex of ${\Gamma_{k+1}}({\cal M})$", is the truncated split
exact differential sequence in the first row of the following
commutative  diagram of split exact sequences:
$$
\begin{array}{rcccccccl}
1&\longrightarrow&Aut({\cal M})&
\mapright{j_{k+1}}&{\Gamma_{k+1}}({\cal M})&
\mapright{\twinup}&{{\cal C}^1_{k+1}}(T{\cal M})&\mapright{\twindp}&
{{\cal C}^2_{k+1}}(T{\cal M})\cr
&&\left|\right|&&\left|\right|&&\mapdown{ad}&&\mapdown{ad}\vespace\cr
1&\longrightarrow&Aut({\cal M})&
\mapright{j_{k+1}}&{\Gamma_{k+1}}({\cal M})&
\mapright{\twinu}&{Der{}^1_{\Sigma,k}}&\mapright{\twind}&
{Der{}^2_{\Sigma,k}}\vespace
\end{array}
$$
where $ad$ is the isomorphism of graded Lie algebras given in the
definition  (\ref{def10}), i.e.
$\forall\tau,\chi~\in~{{\cal C}^{\bullet}_k}(T{\cal M})$:
$$
ad(\grobral\tau,\chi\grobrar')=[ad(\tau),ad(\chi)].
$$\car
\label{def13}
\end{definition}

These sequences can be restricted to ${\widehat{\cal{R}}}_{k+1}$
if ${\widehat{\cal{R}}}_{k+1}$ satisfies the same properties as
those given in the sequel of definition (\ref{def3}) concerning
$D'$, and moreover if it is 2-acyclic.

\addcontentsline{toc}{subsection}{5.2. The gravitational and electromagnetic fields}
\subsection*{\centering 5.2. The gravitational and electromagnetic fields}

\hskip 6mm On the one hand, one has the following commutative diagram:

$$
\begin{array}{ccc}
\Jus&\mapright{\twinsu}&\tsms\otith\Jus\cr
\mapup{\phi_0}&&\mapup{\phi_1}\vespace\cr
\rdc&\mapright{\twinupc}&{\widehat{\cal C}^1_2}\vespace\cr
\mapvegal{id}&&\mapdown{ad}\vespace\cr
\rdc&\mapright{\twinuc}&{Der{}^1_{\Sigma,a}}(\widehat{B}_1)\vespace
\end{array}
$$

where ${\widehat{B}_1}=\gtms{}\otith\Jus$, and on the other hand, one
obtains  the following relation deduced from the relations
(\ref{30}) and (\ref{31}):
$\forall\,{\alpha_1}\in\Jus\equiv{\widehat{B}{}^0_1}$,

$$
\twinsu({\alpha_1})=D'({\alpha_1})-\twinuc(\fdc)({\alpha_1}).
$$

Also let us define $\dgou{1}\equiv\twinuc(\fdc)
\in{Der{}^1_{\Sigma,a}}(\widehat{B}_1)$ and $\mbox{\goth D}$ such
that:

$$
\mbox{\goth D}=D'-\dgou{1}\in{Der{}^1_{\Sigma,a}}(\widehat{B}_1).
$$

It follows one can rewrite:

$$
\mbox{\goth D}\circ\mbox{\goth D}=-\twindc(\dgou{1}),
$$

where $\twindc$ is the $2^{nd}$ non-linear Spencer operator of the
sophisticated Spencer complex of $\rdc$:

$$
1\longrightarrow{Aut({\cal M})}\longrightarrow\rdc\mapright{\twinuc}
{Der{}^1_{\Sigma,a}}(\widehat{B}_1)\mapright{\twindc}
{Der{}^2_{\Sigma,a}}(\widehat{B}_1).
$$

Then, from the definition of $Der{}^{\bullet}_{\Sigma,k}$, there
exists $\tc{2}\in{\widehat{\cal C}^1_2}$ such that
$\dgou{1}=ad(\tc{2})$ and therefore one can write also:

$$
\mbox{\goth D}\circ\mbox{\goth D}=-ad(\twindpc(\tc{2})).
$$

Hence $\forall\,{\alpha_1}\in\Jus$ one deduces the relation:

$$
\mbox{\goth D}
\circ\twinsu({\alpha_1})=-ad(\twindpc(\tc{2}))({\alpha_1}),
$$

and also $\forall\tc{2}\in{\widehat{\cal C}^1_2}$:

$$
\begin{array}{rcl}
\mbox{\goth D}\circ{\phi_1}(\tc{2})&=&D'\circ{\phi_1}(\tc{2})-
\dgou{1}\circ{\phi_1}(\tc{2})\cr
&=&D'\circ{\phi_1}(\tc{2})-ad(\tc{2})\circ{\phi_1}(\tc{2})\cr
&=&D'\circ{1\over{n}}Tr^1(\tc{2})-
ad(\tc{2})\circ{\phi_1}(\tc{2})\cr
&=&{1\over{n}}Tr^1\circ{D'}(\tc{2})-
ad(\tc{2})\circ{\phi_1}(\tc{2})\cr
&\stackrel{def.\,{\phi_2}}{\equiv}&{\phi_2}\circ{D'}(\tc{2})-
ad(\tc{2})\circ{\phi_1}(\tc{2}),
\end{array}
$$

which can be rewritten:

$$
\mbox{\goth D}\circ{\phi_1}(\tc{2})={\phi_2}(\twindpc(\tc{2}))+
{1\over{2}}{\phi_2}(\grobral\tc{2},\tc{2}\grobrar)-
ad(\tc{2})\circ{\phi_1}(\tc{2}),
$$

where ${\phi_2}$ is the $\theta$-linear morphism:

$$
{\phi_2}:{\widehat{\cal C}^2_2}\longrightarrow
{Der{}^2_{\Sigma,1}}\equiv{Der{}^2_{\Sigma,a}}(\widehat{B}_1),
$$

satisfying the relations:
$\forall\,{\displaystyle{\widehat{\sigma}^{(2)}}}
\in{\widehat{\cal C}^2_2}$,

$$
\begin{array}{c}
{\phi_2}({\displaystyle{\widehat{\sigma}^{(2)}}})\equiv
{1\over{n}}Tr^1({\displaystyle{\widehat{\sigma}^{(2)}}})\cr
\mbox{and}\quad{\phi_2}\circ{D'}=D'\circ{\phi_1}.
\end{array}
$$

Then , defining
$\twinsd\stackrel{def.}{\equiv}\mbox{\goth D}
{\Big/_{\tsms\otith\Jus}}$, we deduce the theorem:

{\bf Theorem} {\it The following diagram of differential sequences
is commutative, i.e.
$\twinsd~\circ~{\phi_1}~=~{\phi_2}~\circ~\twindpc$:}

\begin{equation}
\begin{array}{rcccccccl}
&&1&&1&&0&&0\cr
&&\mapdown{}&&\mapdown{}&&\mapdown{}&&\mapdown{}\vespace\cr
1&\longrightarrow&\Gamma_G&\mapright{j_2}&{\cal R}_2&
\mapright{\twinup}&{\cal C}^1_2&\mapright{\twindp}&{\cal C}^2_2\vespace\cr
&&\mapdown{}&&\mapdown{}&&\mapdown{}&&\mapdown{}\vespace\cr
1&\longrightarrow&\Gamma_{\widehat{G}}&\mapright{j_2}&\rdc&
\mapright{\twinupc}&{\widehat{\cal C}}^1_2&\mapright{\twindpc}&
{\widehat{\cal C}}^2_2\vespace\cr
&&\mapdown{}&&\mapdown{\phi_0}&&\mapdown{\phi_1}&&
\mapdown{\phi_2}\vespace\cr
1&\longrightarrow&\theta&\mapright{j_1}&\Jus&
\mapright{\twinsu}&\tsms\otith\Jus&\mapright{\twinsd}&
\gtms{2}\otith\Jus\vespace\cr
&&&&\mapdown{}&&\mapdown{}&&\mapdown{}\vespace\cr
&&&&1&&0&&0\vespace\cr
\end{array}\label{32}
\end{equation}

{\it if the condition below is satisfied:
$\forall\,{\widehat{\tau}}\in{{\widehat{\cal C}}^1_2}$,}

\begin{equation}
{\phi_2}(\grobral{\widehat{\tau}},{\widehat{\tau}}\grobrar)=
2\,ad({\widehat{\tau}})\circ{\phi_1}({\widehat{\tau}}).
\label{33}
\end{equation}\car

Finally, one notices from a diagram chasing in the
diagram (\ref{32}) that the sequence is exact at $\tsms\otith\Jus$
if and only if one restricts to $Im\,\twinsu$, i.e. one has only the
short split exact sequence:

\begin{equation}
1\longrightarrow\theta\mapright{j_1}\Jus
\mapright{\twinsu}\twinsu(\Jus)\mapright{\twinsd}0.
\label{33bis}
\end{equation}

This might be (?) an illustration of the fondamental theorem of Spencer
on the deformations of Lie structures
(Spencer 1962, 1965, Goldschmidt 1976a, 1976b, 1978a, 1978b, 1981)
since $\Jus$ is endowed with a Lie group structure.

Let us denote $\widehat{\cal C}$ and $\pi^{(2)}$ by:

$$\begin{array}{lcl}
{\widehat{\cal C}}&\equiv&(\pve,\ptge)\in
{\widehat{B}}{}^1_1=\tsms\otith\Jus\cr\\
\pi^{(2)}&\equiv&(\Omega^{(2)},{\cal K}^{(2)})\in
{\widehat{B}}{}^{1\,\ast}_1\otith
({{\widehat{\cal C}}^1_2}/{{\cal C}^1_2})\cr\\
\mbox{and}&&<{\widehat{\cal C}}|\pi^{(2)}>
\stackrel{def.}{\equiv}<\pve|\Omega^{(2)}>+
<\ptge|{\cal K}^{(2)}>={\widehat{\tau}}-\tau^{(2)}
\end{array}$$

with, as a consequence of the relations (\ref{30}) and (\ref{31}):

$$
{\phi_1}(\pi^{(2)})=1d_{{\widehat{B}}{}^1_1}\in
{{\widehat{B}}{}^{1\,\ast}_1}\otith{{\widehat{B}}{}^1_1}
\mbox{\quad{and}\quad}{\phi_1}(\tau^{(2)})=0.
$$

Then, from (\ref{33}), one obtains the following equivalent
relations:

\begin{eqnarray}
{\phi_2}(\grobral\tau^{(2)},\tau^{(2)}\grobrar)&=&0\label{34.1}\\
{\phi_2}(\grobral\tau^{(2)},\pi^{(2)})\grobrar)&=&ad(\tau^{(2)})\label{34.2}\\
{\phi_2}(\grobral\pi^{(2)},\pi^{(2)}\grobrar)&=&2\,ad(\pi^{(2)}),\label{34.3}
\end{eqnarray}

which are the defining constraints on the $\tau^{(2)}$ and the
susceptibilities $\pi^{(2)}$.\par
In an orthonormal system of local coordinates, one can write
$\twinsd({\widehat{\cal C}})\equiv(\widehat{\cal G},\widehat{\cal H})=
{\phi_2}({\widehat{\sigma}}^{(2)})\in\gtms{2}\otith\Jus=
{{\widehat{B}}{}^2_1}$ in the form ($h,i,j,k,l,r,s=1,...,n$):

$$\begin{array}{rcl}
\widehat{\cal G}_{[i\,j]}&=&{1\over{n}}(
{\widehat{\sigma}}{}^k_{k,[i\,j]}+
{\widehat{\sigma}}{}^k_{\,,[i\,j]}{\gamma{}^h_{h\,k}}
)\cr\\
&\equiv&{{\widehat{B}}{}^k_{[\,i}}({{\partial}{}_{|\,k\,|}}
{\pve_{j\,]}})-{\widehat{\cal F}_{[i\,j]}}-
{\widehat{\tau}}{}^k_{[j,i]}{\pve_k}\cr\\
\widehat{\cal H}_{j,[k\,i]}&=&
{1\over{n}}\left(
{\widehat{\sigma}}{}^h_{h\,j,[k\,i]}+
{\widehat{\sigma}}{}^l_{j,[k\,i]}{\gamma{}^h_{h\,l}}+
{\widehat{\sigma}}{}^k_{\,,[k\,i]}({{\partial}{}_j}
{\gamma{}^h_{h\,k}})
\right)\cr\\
&\equiv&{{\widehat{B}}{}^h_{[\,k}}({{\partial}{}_{|\,h}}
{\ptge_{j,|\,{i}\,]}})-{\widehat{\cal E}_{j\,[k,i]}}-
{\widehat{\tau}}{}^r_{[i,k]}{\ptge_{j,r}}-
{\widehat{\tau}}{}^r_{j,[\,k}{\ptge_{|\,r,|\,i\,]}}-
{\widehat{\tau}}{}^r_{j\,[i,k]}{\pve_{r}}
\end{array}$$

where

\begin{enumerate}
\item $\widehat{\cal F}\in\gtm{2}$ is the skew-symmetric
part of $\ptge$ one calls the electromagnetic field tensor
(or the  Faraday tensor),
\item $(\widehat{\cal E},\ptge,\pve)\in\tsms\otith
{\displaystyle\not{\!\!{R_{{}_2}}}}$, i.e. $\forall\,X,Y,Z\,\in
T{\cal M}$, $\forall\,c_0\in{\R}$,
$$
<\widehat{\cal E}(X,Y)|{Z}>=
{1\over{2}}<\ptge({\nabla_X}Y+{\nabla_Y}X)|{Z}>-{c_0}\,
\omega(X,Y)<\pve|{Z}>,
$$
\item $\tc{2}$ satisfies the relation (\ref{31bis}).
\end{enumerate}

The terminology one uses for $\widehat{\cal F}$,  comes
obviously from the previous relations analogous to those obtained
in the Maxwell theory. Indeed, since
$(\widehat{\cal G},\widehat{\cal H})=0$ because  of
the exactness condition in the sequence (\ref{33bis}), one deduces:

$$
{\widehat{\cal F}_{[i\,j]}}={{\widehat{B}}{}^k_{[\,i}}
({{\partial}{}_{|\,k\,|}}{\pve_{j\,]}})+
{\widehat{\tau}}{}^k_{[i,j]}{\pve_k},
$$

and from the expression for $\widehat{\cal H}$:

$$
{\widehat{B}{}^h_{[\,k}}({{\partial}{}_{|\,h}}
{\widehat{\cal F}_{|\,j,i\,]}})+
{\widehat{\tau}{}^r_{[k,j}}{\widehat{\cal F}_{i,r]}}=0.
$$

In the case of the weak fields limit, the latter become:

$$\left\{
\begin{array}{rcll}
\widehat{\cal F}&=&\de\pve&\mbox{\quad{(Faraday tensor)}}\cr
\de\widehat{\cal F}&=&0&\mbox{\quad{(Bianchi identity).}}
\end{array}
\right.
$$

On the contrary, the symmetric part
$\widehat{\cal P}\in{{S}_2}\tsm$ of $\ptge$, called ``the
gravitational field tensor", satisfies a first order PDE depending
on $(\pve,\ptge)$, and thus $\widehat{\cal P}$ varies
even if there is only an electromagnetic field! One has to notice
that the symmetric part $({\partial_i}{\pve_j})$ is never taken into
account in physics, in contradistinction to $\pve$. As a result all
its derivatives should be physical observables. Finally, one deduces
from the exactness of the complex (\ref{33bis}) that no current  of
magnetic charges can exist. In other words, the Bianchi identity
must be satisfied. Hence, the lack of magnetic charges can be
justified in the framework of the Spencer cohomology
of conformal Lie structures but not in the de Rham  cohomology
framework.

\addcontentsline{toc}{section}{6. The dual linear
Spencer complex and the Janet complex}
\section*{\centering 6. The dual linear Spencer complex
and the Janet complex}

\addcontentsline{toc}{subsection}{6.1. The dual linear
Spencer complex}
\subsection*{\centering 6.1. The dual linear Spencer complex}

\hskip 6mm One shall refer in this chapter to the definitions given in
references (Goldschmidt 1976a, 1976b, 1978a, 1978b, 1981, Gasqui {\it et al.}  
1984, Pommaret 1995) relating to the dual linear
Spencer complex. Mainly,
in this chapter, one has to build up the dual operator $\ditwinsu$ of
the infinitesimal operator $\itwinsu$ of $\twinsu$, and the morphisms
$L_{n-i}$ ($i=0,1$) in the commutative diagram:

$$
\begin{array}{rcccccccl}
0&\longrightarrow&\theta_{\cal M}&\mapright{j_1}&\Jus&
\mapright{\itwinsu}&\itwinsu(\Jus)\subseteq\cs{1}{1}&
\longrightarrow&0\cr
&&&&\mapdown{L_n}&&\mapdown{L_{n-1}}&&\rule{0pt}{20pt}\cr
&&&&\as{n}{1}&\mapleft{\ditwinsu}&\as{n-1}{1}&\mapleft{}&0
\rule{0pt}{20pt}
\end{array}$$

where the $\as{n-i}{1}$ are the dualizing fiber bundles
($\cs{0}{1}\equiv\Jus$):

$$\as{n-i}{1}=\gtms{n}\otith\dcs{i}{1}.$$

One uses a rather classical method to determine the adjoint operators
on a connected and oriented compact without boundaries (see P. J. Olver
for instance 1986). To this purpose, one gets before a conformally
equivariant Lagrangian density {\goth L}:

$$
\mbox{\goth L}:\Jus{\times_{{}_{\cal M}}}\cs{1}{1}\equiv\cs{0}{1}
{\times_{{}_{\cal M}}}\cs{1}{1}\longrightarrow\gtms{n},
$$

and one defines the morphisms $L_{n-i}$ and $\ditwinsu$ integrating by
parts the infinitesimal variation of {\goth L}, with $\itwinsu$ defined
by (${\alpha_1}\equiv(\alpha,\beta)$):

$$
\left\{
\begin{array}{rcl}
\pve&=&\de\alpha-\beta\cr
\ptge&=&\de\beta-\mu
\end{array}
\right\}\equiv\itwinsu{\alpha_1}
$$

where $\forall\,X,Y\in{{\cal C}^1}(T{\cal M})$:

$$
\mu(X,Y)={1\over{2}}\beta({\nabla_X}Y+{\nabla_Y}X)-{c_0}\,\alpha\,\omega
(X,Y).
$$

Let us define in $\as{n-1}{1}$, the images of $(\pve,\ptge)$ by the
morphism $L_{n-1}$ with the relations:

$$\left\{
\begin{array}{rcl}
\widehat{\cal J}&=&(\partial\mbox{\goth L}/\partial\pve)
\mbox{\quad(electric current)}\cr
\\
\widehat{\cal N}&=&(\partial\mbox{\goth L}/\partial\ptge),
\end{array}
\right.$$

and in $\as{n}{1}$, the images of $\alpha_1$ by the
morphism $L_{n}$ with the relations:

$$\left\{
\begin{array}{rcl}
\widehat{\cal S}&=&(\partial\mbox{\goth L}/\partial\alpha)
\mbox{\quad(density of entropy)}\cr
\\
\widehat{\cal Q}&=&(\partial\mbox{\goth L}/\partial\beta).
\end{array}
\right.$$

One obtains from the infinitesimal variation $\delta\mbox{\goth L}$,
the following defining relations of $\ditwinsu$:

$$\left\{
\begin{array}{rcl}
\widehat{\cal S}&=&\div\widehat{\cal J}-{c_0}\,
<\omega|\widehat{\cal N}>\cr
\\
\widehat{\cal Q}&=&\widehat{\cal J}+{\div_2}\widehat{\cal N}+
<\zeta|\widehat{\cal N}>,
\end{array}
\right.$$

where ${\div_2}$ and $\zeta$ are morphisms such as:

$$\begin{array}{rccl}
{\div_2}:&\stackrel{2}{\otimes}T{\cal M}&\longrightarrow&T{\cal M}\cr
&u\otimes{v}&\longrightarrow&v \div(u) - u \div(v) + [u,v]
\end{array}$$
and
$$\begin{array}{rccl}
\zeta:&\stackrel{2}{\otimes}{{\cal C}^1}(T{\cal M})&
\longrightarrow&{{\cal C}^0}(T{\cal M})\cr
&u\otimes{v}&\longrightarrow&
\zeta(u\otimes{v})=\gamma(u)v
=\gamma(v)u\equiv<\xi\mid u\otimes\,v>.
\end{array}$$

From these results, a further step would be to know what the relations
between $\ditwinsu$ and the thermodynamical laws for
the irreversible processes are, since the Lagrangian density {\goth L}
depends on temperature and variables of internal states.
Hence, {\goth L} might be identified with a Gibbs free enthalpy function,
but moreover, setting certain  conditions, the conformal equivariance
laid down to {\goth L}, might also allow us  to consider {\goth L} as
a wave-function, as we shall see in what follows.

\addcontentsline{toc}{subsection}{6.2. The Janet complex of the
Lagrangian density {\goth L}}
\subsection*{\centering 6.2. The Janet complex of the
Lagrangian density {\goth L}}

\hskip 6mm Again, one refers for the definitions to
Pommaret's previous papers (Pommaret 1989, 1994) about the Janet complex and
also to Gasqui-Goldschmidt's ones (Gasqui {\it et al.} 1988),
but with an ``un-named complex".
The conformal equivariance of {\goth L} will merely allow us to  obtain
explicitly the first Janet operator
${{\D}_1}:{{\cal F}_0}\longrightarrow{{\cal F}_1}$ where
${\cal F}_0$ stands for the line bundle over $\cs{0}{1}
{\times_{{}_{\cal M}}}\cs{1}{1}$:

$$\begin{array}{c}
{\cal F}_0\cr
\mapdown{}\mapup{\mbox{\goth L}}\cr
\cs{0}{1}
{\times_{{}_{\cal M}}}\cs{1}{1}
\end{array}$$

From now on, one presents a dependent coordinates formulation choosing
an orthonormal system of local coordinates.
To make ${\D}_1$, one must find the transformation rule of
{\goth L} by the action of the conformal Lie pseudogroup. Namely, if
$\fdc\in\rdc$ and $\ftci{2}={\displaystyle{{\hat{g}}_2}}$, one has
the transformation rules:

$$\left\{
\begin{array}{rcl}
y&=&\hat{f}(x)\simeq{x}+\hat{\xi}(x)\cr
\\
\alpha'\circ\hat{f}&=&\alpha-{1\over{n}}\ln|{\det{J}(\fuc)}
|\cr
\\
{\beta'_i}\circ\hat{f}&=&{\hat{g}{}^j_i}\circ\hat{f}\,[
{\beta_j}-{1\over{n}}{\hat{g}{}^k_l}\circ
\hat{f}\,{\hat{f}{}^l_{k\,j}}]\cr
\\
{\pve'_i}\circ\hat{f}&=&{\hat{g}{}^j_i}\circ\hat{f}\,{\pve_j}\cr
\\
{\ptge'_{i,j}}\circ\hat{f}&=&{\hat{g}{}^r_i}\circ\hat{f}\,\,
{\hat{g}{}^s_j}\circ\hat{f}\,{\ptge_{r,s}}
\end{array}
\right.$$

and one must have:

\begin{equation}
\mbox{\goth L}(\hat{f},\alpha'\circ\hat{f},
\beta'\circ\hat{f},\pve'\circ\hat{f},\ptge'\circ\hat{f})\,\det
{J}(\fuc)=\mbox{\goth L}(x,\alpha,\beta,\pve,\ptge).
\label{el}
\end{equation}

One obtains the infinitesimal condition and the definition of
${\D}_1$:

$$
{{\D}_1}(\mbox{\goth L})\equiv\left\{
v^\mu(\mbox{\goth L})+(\div\,\xi^\mu)\,\mbox{\goth L}\equiv
\mbox{\goth K}\right\}=0,
$$

whatever is $v^\mu$ ($\mu=1,...,\dim\mbox{\goth g}_c$), generator of the
Lie algebra $\mbox{\goth g}_c$ of the conformal Lie group:

$$
v^\mu={\hat{\xi}^{\mu,j}}{\partial_j}+
{\phi{}^\mu}(\alpha){\partial_\alpha}+
{\phi{}^\mu}(\beta_j){\partial_{\beta_j}}+
{\phi{}^\mu}(\pve_j){\partial_{\pve_j}}+
{\phi{}^\mu}(\ptge_{k,l}){\partial_{\ptge_{k,l}}},
$$

where

$$\left\{
\begin{array}{rcl}
{\phi{}^\mu}(\alpha)&=&-{1\over{n}}{\hat{\xi}^{\mu,k}_k}\cr
\\
{\phi{}^\mu}(\beta_j)&=&-[{\hat{\xi}^{\mu,k}_j}{\beta_k}+
{\hat{\xi}^{\mu,k}_{k\,j}}]\cr
\\
{\phi{}^\mu}(\pve_j)&=&-{\hat{\xi}^{\mu,k}_j}{\pve_k}\cr
\\
{\phi{}^\mu}(\ptge_{k,l})&=&-[{\hat{\xi}^{\mu,h}_k}{\ptge_{h,l}}+
{\hat{\xi}^{\mu,h}_l}{\ptge_{k,h}}]
\end{array}
\right.$$

and ${\hat{\xi}^\mu_2}\in{\widehat{R}_2}$.
Moreover, $\mbox{\goth K}\equiv(\mbox{\goth K}^\mu)$ is transformed like
{\goth L} and therefore the zero section of ${{\cal F}_1}$ is defined
with $\dim{\mbox{\goth g}_c}=\dim{{\cal F}_1}$ constants $c^\mu$ such that
$\mbox{\goth K}^\mu=c^\mu\mbox{\goth L}$ meaning that (\ref{el}) (or the
Lagrangian density) is defined
up to a multiplicative constant. The second Janet operator
${\D}_2$ such as ${\D}_2\circ{{\D}_1}=0$
and ${{\D}_2}:{{\cal F}_1}\longrightarrow{{\cal F}_2}$, is
defined from the involution of the generators of $\mbox{\goth g}_c$:

$$
[v^\mu,v^\nu]={c^{\mu\,\nu}_{\quad\lambda}}\,v^\lambda,
$$

where the ${c^{\mu\,\nu}_{\quad\lambda}}$ are the constants of
structure of the $\mbox{\goth g}_c$ (see P. J. Olver 1986 for a precise
definition of the bracket used above).
Thus, one obtains:

$$
{{\D}_2}(\mbox{\goth K})\equiv\left\{
v^\mu(\mbox{\goth K}^\nu)-
v^\nu(\mbox{\goth K}^\mu)-
{c^{\mu\,\nu}_{\quad\lambda}}\,\mbox{\goth K}^\lambda+
(\div\,{\hat{\xi}^\mu})\,\mbox{\goth K}^\nu-
(\div\,{\hat{\xi}^\nu})\,\mbox{\goth K}^\mu
\right\}.
$$

Leading down $\mbox{\goth K}^\mu=c^\mu\mbox{\goth L}$ in the latter
expression, one deduces the following constraints on the constants
${c^{\mu\,\nu}_{\quad\lambda}}$ of structure:

$${c^{\mu\,\nu}_{\quad\lambda}}{c^\lambda}=0,$$

with $\displaystyle{{\D}_1}(\mbox{\goth L})^\mu=
{c^\mu}\mbox{\goth L}.$

This latter equation gives in fact only one arbitrary constant among
the ${c^\mu}$'s and perhaps might be ascribed to a constant such
as the Planck one.

As a particular case, if the metric $\omega$ and $\pve$ are constants,
$\alpha$ and $\beta$ vanish and {\goth L} function of $x$ and $\pve$
only, such that there exists $\eta$ satifying the relation

$$
\widehat{\cal J}^k=\eta^k\mbox{\goth L},
$$

then $\forall\,\mu,h$ one has ${\hat{\xi}^{\mu,k}_k}=
{\hat{\xi}^{\mu,k}_{k\,h}}=0$ and

$$
\left[
{\hat{\xi}^{\mu,k}}{\partial_k}-(\eta^k{\hat{\xi}^{\mu,h}_k})
{\pve_h}
\right]\mbox{\goth L}-c^\mu\,\mbox{\goth L}=0.
$$

Thus, one obtains an analogous PDE to the Dirac equation but being
equivariant and not only covariant like the Dirac equation.
Nevertheless {\goth L} being a real function, it can be interpreted
as  a wave-function only if one has a definition of a measure of
probability. We propose two ways of doing this in the conclusion,
and we make suggestions for taking the weak and
strong interactions into account in the model.

\addcontentsline{toc}{section}{7. Conclusion}
\section*{\centering 7. Conclusion}

\hskip 6mm The function {\goth L} being real, let us suggest a notion
of mesure of probability. Two of them can be proposed. The first one
is related to an
``\`{a} la Misra-Prigogine-Courbage" (MPC) (Misra {\it et al.} 1979, Misra  
1987) approach. With more details, the
wave-function {\goth S} of the model would be the function obtained when
applying the integral action of the evolution operator to
an initial Lagrangian density {\goth L}. First, we thus obtain a wave
interpretation. Second, within the MPC theory, the reduction of
the wave-paquet (or the collapse) would axiomatically be defined as
the achievement of a K-partition of  functions {\goth S}
on a compact (on which the initial  Cauchy data would be defined)
of a submanifold of the Minkowski  space-time to determine a single
function {\goth S}.  This K-partition of functions {\goth
S} would then be a set of new initial Lagrangian  densities
{\goth L}  before any time evolution. The physical  measurement
thus achieved  would not lead us to obtain pure states and with the
assumption (to be confirmed) that the conformally equivariant
Lagrangian densities {\goth L} would be Kolmogorov flows, it would
then be possible, according to the MPC theory, to build up a
non-commutative algebra of observables on the analogy of quantum
mechanics. Moreover, the K-partition would be achieved on a set of
Lagrangian densities of type (m,s) deduced by projection on the
basis of the eigenstates of the Poincar\'{e} group during the
process of measurement (from a certain point of view, the apparatus
of measurement, being themselves of the type (m,s)). Thus, as a
result of this projection, we shall see a particle state and
initial Lagrangian densities of type (m,s) not depending on the
fields of interaction anymore. The particle states might be
interpreted as a split up of variables between the variable of
position $x$ and the fields variables (during specific physical
measurements). In other words, one could say that the energy of the
fields of interaction would be transfered by radiation to the
apparatus of measurement. From a more general point of view, the
process of ``fragmentation" would be defined by the change from a
tensorial Hilbert space to a decomposable Hilbert space (${\cal
H}_i$; $i=1,2$: Hilbert spaces):

$$
{\cal H}_1\otimes{\cal H}_2\longrightarrow{\cal H}_1\oplus
{\cal H}_2
$$

Process that might be interpreted as well as a separation of phase in
thermodynamics for the function {\goth L} is  similar to a Gibbs
function. Under these conditions, the analytical developments of
the initial Lagrangian densities define  tensors of
susceptibility directly associated with the ``classical" states of the
condensed matter, therefore with the classical notion of macroscopic
particles. Thus, the macroscopicity would be associated with the
``degree" of separability of the initial densities {\goth L}.
From our point of view, a second possibility used to
define a measure of probability would be to refer to the Penrose
approach (Penrose {\it et al.} 1986) on the complexification of the wave-function
in the frame of the theory of twistors, with moreover, the notion
of measure of probability defined as it usually is in the first
quantization. This procedure is completly different from the one
consisting in simply making an algebraic extension from the field of
reals to the field of  complexes, of the  wave-functions satisfying
PDE on  the field of  reals.
Such an extension would lead to complex solutions whose real
and imaginary parts would each separatly be solutions.
Consequently, we assume that a complexification presenting a real
interest would give access to the determination of complex solutions
which, on the contrary, would be neither their real parts, nor
their imaginary ones, solutions for the given system of PDE.  To
achieve that, Penrose suggested decomposing a problem defined on
the real Minkowski space ${\cal M}$ into a set of physical subsystems
defined on submanifolds of ${\cal M}$ and having at least one
(complex) spin structure, like the celestial spheres $S^{\pm}$ are,
or  at fixed time, the sphere $S^{3}$. The spin structures are
associated with vector bundles like the  tangent spaces, or on a
rather similar manner  (up to the target manifold), the 1-jet
bundle ${J_1}(\epsilon)$ for example. In the present situation, we
shall refer to the $\itwinsu(\Jus)\subseteq\cs{1}{1}$ bundle over $\Jus$.
We therefore must present the tensors $\pve$ and $\ptge$ with spinors.
In the more  simple case with only one spin structure (on $S^{2}$
for example) and taking up again the Penrose  indexed notation, one
need at least 3 spinors to decompose the  tensor $\pve$:

$$
{\pve_a}={\alpha_A}.{\overline{\alpha}_{A'}}+
{\beta_A}.{\overline{\beta}_{A'}}-
{\gamma_A}.{\overline{\gamma}_{A'}}.
$$

For only under these conditions, $\pve$ is of any norm, but with a
spinor, $\pve$ must be of the type light and with two spinors, the
norm of $\pve$ must be non-negative or non-positive.  As far as
I know this simple decomposition has  never
been presented yet (apart from the one with 2 spinors or 3
maximum and with a + mark). But it has indeed a few avantages.
First, it is associated (Bars {\it et al.} 1990) with the representations  
$\{3\}$ and $\{\bar{3}\}$ of $SU(2,1)$ isomorphic to the one of $SU(3)$. It is
only about representations; the group $SU(2,1)$ not being a  dynamical
gauge group of the model, like for example the Yang-Mills type, but only
a group of internal classification and of invariance of the
decomposition. The non-compactness of the group cannot appear
determining to us within this context. But its finished unitary
irreducible representations constructed from the two fundamental
representations $\{3\}$ and $\{\bar{3}\}$ only do. Moreover, a
symmetry breaking from $SU(2,1)$ to $SU(2){\times}U(1)$ or
$SU(1,1){\times}U(1)$ is associated to a symmetry breaking of type
$T$ (for instance during the process of fragmentation) and
$\pve$ is then of a non-negative or non-positive norm. All of that
suggests a possible link with a theory of weak and strong
interactions. Before that, formally deriving the spinors
decomposing $\pve$, we obtain the following  decomposition for
$\ptge$ (Penrose {\it et al.} 1986, see chapter 4.4.7.) using the Leibniz law
for the  spinors:

$$
{\ptge_{a,b}}={\alpha_A}.{\overline{\Gamma}_{B'A',B}}+
{\overline{\alpha}_{A'}}.{\Gamma_{BA,B'}}+
{\beta_A}.{\overline{\Theta}_{B'A',B}}+
{\overline{\beta}_{A'}}.{\Theta_{BA,B'}}+
{\gamma_A}.{\overline{\Omega}_{B'A',B}}+
{\overline{\gamma}_{A'}}.{\Omega_{BA,B'}},
$$

where $\Gamma$, $\Theta$ and $\Omega$ are an any mixed spinors of
valence (2,1). They can decompose themselves into irreducible
spinors with respect to $SL(2,{\C})$:

$$
{\Gamma_{AB,B'}}={\Gamma_{(AB),B'}}-{1\over{2}}
{\epsilon_{[AB]}}{\Gamma_{C\,\,\,,B'}^{\,\,\,\,C}}
$$

(with similar expressions for  $\Theta$ and $\Omega$). Let us
notice  that at a given fixed time $t$ and with ($\alpha$,$\beta$) fixed
(i.e. during the measure), $SL(2,{\C})$ being associated with
${O{}^{{}^{{\hskip -.3mm}+}}_\uparrow}(1,3)$, the irreducibility is
thus, according to $SU(2)$, locally with respect to $O(3)$. Moreover,
${\Gamma_{(AB),B'}}$ is defined by 6 complex components, linearly
independant and ${\Gamma_{C\,\,\,,B'}^{\,\,\,\,C}}$ by 2  complex
components. The interest of such a decomposition appears as soon
as the assumption is made that the Lagrangian density {\goth L} is
holomorphic in the various spinors on the submanifolds endowed with
spin structures. Then {\goth L} satisfies  complex PDE deduced from
the real PDE defined on ${\cal M}$ by lifting on the spinors
bundle.  But, the physical meaning of the holomorphy is that the
physical  system can precisely be confined on the submanifolds of
{\cal M}  endowed with a spin structure. If we make the physical
interpretation that to the spinors decomposing $\pve$ and $\ptge$
are associated  some fields of interactions or particles, it means
that the latter can only be considered as ``free" if they are
precisely confined on  those submanifolds; the density {\goth L}
being then a complex  solution on the spinors bundle of the various
lifted PDE of the  model, from which neither the real part, nor the
imaginary part are  any solutions.  On that subject, one cannot
help thinking of the very controversial - and  to my knowledge
still un-confirmed -  Larue
{\it et al.} (1977) and Schaad {\it et al.} (1981) experiences  on ``free"  
quarks confined on
bidimensional  ($\simeq{S^{2}}$) layers of superconducting (!)
nobium  covering  microballs of tungsten (see also Goldman 1995, 1996).
Especially because the  decomposition presented above could
suggest that the  simple spinors $\alpha$, $\beta$ and $\gamma$
would be associated  with some fields of fermion of mono-colored
gluons of spin $1/2$  (and not to bi-colored bosons!) and that the
symmetric mixed spinors of valence (2,1) would be 3 quarks states
of spin
$1/2$ determined by the 3 contracted spinors
${\Gamma_{C\,\,\,,B'}^{\,\,\,\,C}}$,
${\Theta_{C\,\,\,,B'}^{\,\,\,\,C}}$ and
${\Omega_{C\,\,\,,B'}^{\,\,\,\,C}}$ on which $SU(2,1)$ acts. At
last, to finish as well as to come back to an application of the
physics of the superconducting states (that some of the theoretical
physicists like Mendelstam (1982) or t'Hooft (1978)
used to explain the quarks confinement) the question is under
which conditions the 4-vector current:

$$\begin{array}{rcl}
\widehat{\cal J}&=&(\partial\mbox{\goth L}/\partial\pve)
\end{array}$$

would be anomalous as well as the Faraday tensor
$\widehat{\cal F}$. A possibility we suggest would be to consider
a third metric $\lambda$ deduced from $\nu$ like $\nu$ is from $\omega$.
The latter metric $\lambda$ would appear in matter such as in crystals
or amorphous materials for instance. Indeed in this case, the second metric
$\nu$ (and its corresponding potentials and fields) would be associated to the
``crystal field" and to a new kind of substrat space-time and new specific
succeptibilities. But in the contrary to $\omega$ the Weyl curvature
associated to $\nu$ is no longer necessarly vanishing as we pointed out in
a previous chapter. Thus, the corresponding Lie equations for conformal
transformations won't be involutive. As a result, the corresponding Spencer
sequences and the relative one won't be exact any more as well.
In particularly magnetic charges might occur and so anyons (Wilczek 1990).
In the framework of symplectic cohomologies of the Lagrangian density
{\goth L}, it can be shown that the anomalies classification is associated to
a second cohomology space knowledge as well as a third one
(Cari\~nena {\it et al.} 1988). The question  remaining to work out would be  
to know if one
can obtain integer cohomologies  from these latter in order to exhibit a
kind of ``quantum structure". To conclude on a more philosophical note,  the
wave-function model we present here could be interpreted  within the
framework of a non-Lavoisian chemistry like G. Bachelard (1940)
developped. In actual fact, the quantons would not be of ``substance",  but
as Bachelard named them, they would rather be ``grains de  r\'{e}action" and
thus, from a certain point of view, like quantons of reactional synthesis or
of separation of phase associated with the concept of fragmentation.  We also
think that this model can't describe global evolution such as  a big-bang
model but only a local evolution. In some ways as  G. Deleuze {\it et al.}  
(1980) say about their striated space (like a substrat space-time)  and smooth  
space concepts
(see chapters: ``12.1227. Trait\'{e} de  nomadologie: la machine de
guerre" and  ``11.1837. De la ritournelle"), ``one can only  know
the path by  exploring it". Also, the form concept for matter in space-time
could  be related to a boundary being the geometrical set of places onto
which the ``grains de r\'{e}action" occur.
\par
\vskip 1cm
{\bf Acknowledgements} I wish to thank Florence Pittolo and Terence Blake who 
very kindly helped me in the translation of the paper and for their helpful
comments and encouragements.

\addcontentsline{toc}{section}{References}
\section*{\centering References}

\vskip 1cm

\begin{verse}
Asher, E. 1973 In  {\it Magnetoelectric Interaction
Phenomena in Crystals.} Proceedings,  Seattle, 69.

Bachelard, G. 1940 {\it La philosophie du non.} Paris: Presse Universitaire  
de France.

Bacry, H. 1967
{\it Le\c{c}ons sur la th\'{e}orie des groupes et les sym\'{e}tries
des particules \'{e}l\'{e}mentaires.} Paris: Gordon and Breach.

Bars, I. \& Zhong-jian, T. 1990
{\it The unitary irreducible representations of SU(2,1).
J. Math. Phys.} (7){\bf 31}, 1576-1586.

Bateman, H. 1910 {\it Proc. Lond. math. Soc.} {\bf 8}, 223.

Cari{\~n}ena, J. F. \& Ibort, L. A. 1988 {\it J. Math. Phys.} {\bf 29}, 541.

Cunningham, E. 1910 {\it Proc. Lond. math. Soc. }{\bf 8}, 77.

Deleuze, G. \& Guattari, F. 1980 {\it Mille Plateaux.} Paris: Les  
\'{e}ditions de minuit.

Dieudonn\'e, J. 1971 {\it \'El\'ements d'analyse, Tome IV.} Paris:  
Gauthier-Villars \'Editeur.

Gasqui, J. \& Goldschmidt, H. 1988
{\it Complexes of Differential Operators and Symmetric Spaces.}
In  {\it Deformation Theory of Algebras and Structures and Applications.} 797-827 
(eds. M. Hazewinkel \& M. Gerstanhaber) Dordrecht-Boston-London: Kluwer  
Academic Publishers.

Gasqui, J. \& Goldschmidt, H. 1984
{\it Deformations infinitesimales des structures conformes plates.}
{\it Progress in Mathematics}, Vol. {\bf 52},  Boston-Basel-Stuttgart: Birkhaeuser.

Goldman, V. J. 1996  {\it Surf. Sci.} {\bf 361}, 1.

Goldman, V. J. \& Su, B. 1995 {\it Science} {\bf 267}, 1010.

Goldschmidt, H. \& Spencer, D. 1976a
{\it On the non-linear cohomology of Lie equations. Part I, Acta Math.} {\bf  
136}, 103-170.

Goldschmidt, H. \& Spencer, D. 1976b
{\it On the non-linear cohomology of Lie equations. Part  II, Acta math.}  
{\bf 136}, 171-239.

Goldschmidt, H. \& Spencer, D. 1978a
{\it On the non-linear cohomology of Lie equations. Part  III, J. Differ.  
Geom.} {\bf 13}, 409-453.

Goldschmidt, H. \& Spencer, D. 1978b
{\it On the non-linear cohomology of Lie equations. Part IV, J. Differ.  
Geom.} {\bf 13}, 455-526.

Goldschmidt, H. 1981
{\it On the nonlinear cohomology of Lie equations. Part V, J. Differ. Geom.}  
{\bf 16}, 595-674.

't Hooft, G. 1978 {\it Nucl. Phys.} B{\bf 138}, 1.

Itzykson, C. \& Zuber, J.-B. 1980
{\it Quantum Field Theory.} McGraw-Hill.

Janner, A. \& Asher, E. 1969 {\it Phys. Lett.} (4)A{\bf 30}, 223.

Janner, A. \& Asher, E. 1978 {\it Physica} {\bf 48}, 425.

Katanaev, M. O. \& Volovich, I. V. 1992 {\it Ann. Phys.} {\bf 216}, 1-28.

Kleinert, H. 1989 {\it Gauge fields in condensed matter.} Vol. {\bf II}.  
World Scientific.

Kumpera, A. \& Spencer, D. C. 1972 {\it Lie Equations, Vol. I: General Theory.
Annals of Mathematics Studies}, Vol. {\bf 73}.

Larue, G. S.\& Fairbank, W. M. \& Hebard, A. F. 1977
{\it Evidence for the Existence of Fractional Charge on Matter*.
Phys. Rev. Lett.} (18){\bf 38}, 1011-1014.

L\'{e}vy-Leblond, J.-M. 1990
{\it ENIGMA OF THE SP$\hbar$IN$x$.}
In {\it Symposium on foundations of modern Physic.} 226-241. World Scientific.

Lur\c{c}at, F. 1964
{\it Quantum Field Theory And The Dynamical Role Of Spin.
Physics} (2){\bf 1}, 95-106.

Mandelstam, S. 1982
{\it Monopoles in Quantum FIeld Theory.}
In {\it Proceedings of the Monopole Meeting.} Trieste, Italy, December 1981.
World Scientific.

Misra, B. \& Prigogine, I. \& Courbage, M. 1979
{\it From deterministic dynamics to probabilistic descriptions.
Physica} A{\bf 98}, 1-26.

Misra, B. 1987
{\it Fields as Kolmogorov Flows. J. Stat. Phys.} (5/6){\bf 48}, 1295-1320.

Olver, P. J. 1986
{\it Applications of Lie Groups to Differential Equations.
Graduate texts in mathematics}, Vol. {\bf 107}. Springer.

Penrose, R. \& Rindler, W. 1986
{\it Spinors and Space-Time.
Cambridge Monographs On Mathematical Physics}, Vol. {\bf I} \& {\bf II}.
Cambridge University Press.

Pommaret, J.-F. 1994
{\it Partial Differential Equations and Group Theory, New Perspectives for  
Applications.
Mathematics and Its Applications}, vol. {\bf 293}.
Dordrecht-Boston-London: Kluwer Academic Publishers.

Pommaret, J.-F. 1995
{\it Suites Differentielles et Calcul Variationnel.
C. R. Acad. Sci. Paris} {\bf 320}, S\'{e}rie I, 207-212.

Pommaret, J.-F. 1989
{\it Gauge Theory And General Relativity.
Reports in Mathematical Physics} (3){\bf 27}, 313-344.

Pommaret, J.-F. 1988
{\it Thermodynamique et th\'{e}orie des groupes.
C. R. Acad. Sci. Paris} {\bf 307}, S\'{e}rie I, 839-842.

Rubin, J. L. 1994
{\it Relativistic Crystalline Symmetry Breaking And Anyonic States in  
Magnetoelectric Superconductors.
Ferroelectrics} {\bf 161}, 335.

Rubin, J. L. 1993
{\it Broken Relativistic Symmetry Groups, Toroidal Moments and
Superconductivity In Magnetoelectric Crystals.
Il Nuovo Cimento D-Cond. Matt. At.} {\bf 15}, 59.

Schaad, L. J. \& Hess Jr., B. A. \& Wikswo Jr., J. P. \& Fairbank, W. M. 1981
{\it Quark chemistry. Phys. Rev.} (4)A{\bf 23}, 1600-1607.

Shih, W. 1986
{\it Une m\'{e}thode \'{e}l\'{e}mentaire pour l\'{e}tude des \'{e}quations aux  
d\'{e}riv\'{e}es partielles.} In
{\it Diagrammes}, {\bf 16}. Paris.

Shih, W. 1987
{\it C. R. Acad. Sci. Paris} {\bf 304}, S\'{e}rie I, 103-106, 187-190, 535-538.

Shih, W. 1991
{\it Gazette. Soci\'{e}t\'{e} Math\'{e}matique de France}, 1.

Spencer, D. C. 1962
{\it Deformation of Structures on Manifolds defined by transitive, continuous 
pseudogroups. Parts I \& II, Annals of Math.} {\bf 76}, 306-445.

Spencer, D. C. 1965
{\it Deformation of Structures on Manifolds defined by transitive, continuous 
pseudogroups. Part III, Annals of Math.} {\bf 81}, 289-450.

Wilczek F. 1990 {\it Fractional Statistics and Anyons Superconductivity.}
World Scientific Publishing Co. Pte. Ltd.

\end{verse}

\end{document}